\documentclass[11pt,twoside,draft,reqno]{birkart}

\usepackage{amssymb}

\usepackage{amsmath}

\usepackage[all,dvips]{xy} 

\newfont{\calligX}{callig15 scaled 666}
\newfont{\script}{callig15 scaled 666}
\newcommand{\nineit}{\small\it}
\newcommand{\ninebf}{\small\bf}
\newcommand{\ninerm}{\small}

\newcommand{\Anu}{A^{\raise-2pt\hbox{$\scriptstyle (\nu)$}}} 
\newcommand{\erfc}{{\rm Erfc}\,}
\newcommand{\re}{{\rm Re}\,}
\newcommand{\im}{{\rm Im}\,}
\newcommand{\Laplacetrafo}{{\text{\script{L\hspace*{4.85pt}}}}}
\newcommand{\la}{{\text{\script{L\hspace*{4.85pt}}}}}
\newcommand{\scriptR}{{\text{\script{R\hspace*{2.85pt}}}}}
\newcommand{\Fgeschwungen}{{\text{\script{F\hspace*{2.85pt}}}}}
\newcommand{\fa}{{\text{\script{F\hspace*{2.85pt}}}}}

\newcommand{\ctrig}{ C_{\rm trig}}
\newcommand{\chyp}[1]{ C_{{\rm hyp}, #1 }}

\newcommand{\grosserAbsatz}{\bigskip\par\noindent}
\newcommand{\mittlererAbsatz}{\medskip\par\noindent}
\newcommand{\kleinerAbsatz}{\smallskip\par\noindent}
\begin{document}
\setcounter{page}{1}

\title[Bessel Processes and Asian options]%
{Bessel processes, the integral of geometric\hfil\break Brownian motion, 
and Asian options}

\author[P. Carr and M. Schr\"oder]{Peter Carr and Michael Schr\"oder}

\address{
Peter Carr\br
Analytic Research at Bloomberg (ARB), 499 Park Avenue, New York, NY 10022\br 
and Courant Institute, NYU, 251 Mercer Street,  New York NY 10012, USA.\br
\br
{\it E-mail address\/}: {\tt pcarr@nyc.rr.com}\br
\br 
Michael Schr\"oder\br 
Keplerstra\ss e 30, D-69469 Weinheim, Germany.\br 
}

\email{schroeder@math.uni-mannheim.de}

\begin{abstract}
This paper is motivated by questions about averages of stochastic 
processes which originate in mathematical finance, originally in 
connection with valuing the so-called Asian options. Starting with
\cite{Yor 92}, these questions about exponential functionals of 
Brownian motion have been studied in terms of Bessel processes 
using the Hartman-Watson theory of \cite{Yor 80}. Consequences 
of this approach for  valuing Asian options proper have been 
spelled out in \cite{Geman+Yor} whose Laplace transform results 
were in fact regarded as a noted advance. Unfortunately, a number 
of difficulties with the key results of this last paper have 
surfaced which are now addressed in this paper. One of them in 
particular is of a principal 
nature and originates with the Hartman-Watson approach itself:
this approach is in general applicable without modifications only 
if it does not involve Bessel processes of negative indices.      
The main mathematical contribution  of this paper is the 
developement of three principal ways to overcome these 
restrictions, in particular by merging stochastics and complex 
analysis in what seems a novel way, and the discussion of their 
consequences for the valuation of Asian options proper. 

\grosserAbsatz
\noindent
{\ninebf Keywords.} Asian options,
                   integral of geometric Brownian motion, Bessel processes,
                   Laplace transform, complex analytic methods in 
                   stochastics.
\end{abstract}
\maketitle
\noindent
{\large\bf 1.\quad Introduction:}\quad 
This paper addresses questions about exponential functionals of 
Brownian motion and the integral of geometric Brownian motion in 
particular. These questions reduce to the study of the quadratic 
variation processes $\Anu$ of geometric Brownian motion which 
for any real drift $\nu$ are explicitly given by the integrals 
over time 
$$
A^{(\nu)}_t= \int_0^t e^{2(\nu w+B_w)}\, dw\, , \qquad
t\in[0,\infty)
$$
with $B$ a standard Brownian motion.  
These processes have both a surprisingly rich theory and manifold 
applications ranging from the physics of random media to 
mathematical finance and insurance. In fact, the insurance motivated 
study of certain perpetuities in \cite{Dufresne90} seems to have 
initiated this line of research. Here, the above integrals over the 
whole time axis are considered and shown to be distributed as the 
reciprocals of certain gamma variables. 
Drawing on his probabilistic interpretation of the Hartman-Watson 
identities in \cite{Yor 80}, Yor was able to extend this work and to
determine the law of the processes $\Anu$ in \cite{Yor 92}.  
This approach, using the Laplace transform and based on  Bessel 
processes, has been found by Yor to open many surprising vistas
on the processes $\Anu$ and their applications, see in particular 
\cite{Yor01}. The interest in these processes $\Anu$ is partially
due to their importance in mathematical finance, in particular 
for understanding the so-called Asian options. 
\mittlererAbsatz
Asian options are widely used financial derivatives.  As discussed 
in Part~I of the paper, these options provide in general nonlinear 
payoffs on the arithmetic average of the price of an underlying asset.  
A common objective in their valuation is to derive an explicit 
expression for a certain functional of $\Anu$. The pursuit of this 
objective has evolved over the last fifteen years as an interplay 
between theoretical and computational perspectives, see \cite{Rogers+Shi} 
or \cite{Dufresne00} for instance.
Yor's work in \cite{Yor 92} clarified the structure of the 
Black-Scholes prices of Asian options by expressing them as certain 
triple integrals. In contrast to this result, however, it was the 
Laplace transform approach of \cite{Geman+Yor}, also based on 
Hartman-Watson theory, which had far reaching consequences for 
the way Asian option valuation is seen today, and on which we focus 
here from Part~II onwards. This is essentially because a very 
explicit expression in terms of Kummer's confluent hypergeometric 
function resulted from this approach for what has been regarded as 
the Laplace transform of the value of such options. From a numerical
point of view, this expression has proved to be amenable to computation 
by numerical inversion, see \cite{Fu+Madan+Wang} as a recent example.
And it seems fair to say that this has lead  financial mathematics 
to a new interest in developing and applying these techniques.  
\mittlererAbsatz
Unfortunately, some difficulties with this Laplace transform 
approach to valuing Asian options have emerged. First, it turned 
out that the Laplace transforms computed in \cite{Geman+Yor} 
are not those of the Asian option's value. This appeared 
to empty the relevance of this result for the finance application 
proper. Luckily, it turned out that there is a reduction of the 
original problem of valuing Asian options to the one considered in 
\cite{Geman+Yor}. All of this is discussed in Part~II of the paper.
\mittlererAbsatz
A difficulty of a more serious and more principal nature, however, 
originates with the Hartman-Watson approach on which the Laplace 
transform computations of \cite{Geman+Yor} are based. Its idea is 
to analyze $\Anu$ using Bessel processes of indices $\nu$, 
and its applicability has limits if $\nu$ is negative because of 
the pathologies which Bessel processes of such negative indices develop. 
Together with background material about Bessel processes,
we have thus given in Part~III a new exposition of  the analysis 
in \cite{Geman+Yor}, which expli\-cit\-ly takes care of the nonnegativity 
restriction on the index $\nu$. Hereby, we have 
been encouraged by the kind support of Yor, and we have tried 
to incorporate his kind tutorials and suggestions. In terms of
financial mathematics, our extension to negative $\nu$ in fact
extends the analysis of Asian options from the zero dividend 
situation originally considered in \cite{Geman+Yor} to one of 
general real dividend yields, and thus to one with general real 
risk-neutral drifts.
In this setting, the condition that $\nu$ is not negative transcribes 
into the postulate that risk-neutral drift is not less than half 
the squared volatility.  Unfortunately, this lower bound on the 
drift restricts the financial applicability of the results. In fact, 
the greater is the volatility, the greater is the range of 
parameters in which the nonnegativity condition on $\nu$ is 
violated and in which the approach does not give the Laplace transforms. 
Unfortunately, it is precisely due to high volatility of the 
underlying asset that Asian options are used in the first place. 
\mittlererAbsatz
Thus, the third contribution of our paper is to remove these 
restrictions by showing the existence of the desired Laplace 
transforms. In fact, we develop three ways for removing these 
restrictions. 
More precisely, we develop three princpal ways of coping with 
the difficulties caused by Bessel processes of negative indices 
in the Hartman-Watson approach. This can be seen as the 
main mathematical contribution of the paper. The approaches 
of Parts~IV and VI are based on an analysis of Bessel processes 
in the spirit of Yor. Our first approach extends the \cite{Geman+Yor} 
analysis to the missing cases using deeper properties of Bessel 
processes of negative indices. Our second approach gives a new 
uniform proof using zero index Bessel processes only. Our
third approach of Part~V extends that of \cite{carr-sch} and 
merges stochastic with complex analytic techniques. This appears 
to be a rather novel and promising line of attack as we are able 
to largely dispense with Bessel processes and focus instead 
on Brownian motion.
\mittlererAbsatz
Apart from all this, our extension of the Laplace transform approach 
of \cite{Geman+Yor} has made possible advances in valuing Asian 
options, some of which are sketched in Part~VII. Hence, it may be 
fair to say that the Laplace transform approach of \cite{Geman+Yor} 
has proved to be a rich source for new results and insights in 
both finance and mathematics. 
\vskip.8cm\noindent
\centerline{\Large\bf I\qquad Prologue}
\vskip.3cm\noindent
%
%
{\large\bf 2.\quad Black-Scholes modelling:}\quad
The results to be discussed originate from the so-called 
risk-neutral approach to the valuation of contingent claims.  
General equilibrium treatments for this and other notions 
developed for the analysis of financial markets and instruments 
are in \cite{Duffie88}, \cite{Duffie92}, \cite[Chapters 1--4]{KSb},  
for instance.
This analysis is based on models of security markets, and this 
section aims to recall the most fundamental of these, 
the Black-Scholes model of security markets. 
\mittlererAbsatz
In fact, we need only that particular case of the Black-Scholes 
model where there are only two securities, and the  understanding 
is that these are traded on markets where their prices are 
determined by equating demand and supply.
First, there is a  riskless security, a bond, whose price $\beta$ grows 
at the continuously compounding positive interest rate $r$, i.e., 
for which we have we have $\beta_t=\exp(rt)$ at any time $t\in [0,\infty)$. 
Then, there is a risky security. The fundamental idea is that all
uncertainties affecting its price $S$ yield a certain probability space. 
In fact, consider for this a complete probability space equipped 
with the standard filtration of a standard Brownian motion  on the 
time set $[0,\infty)$. Giving expression to the fact that $S$ comes 
as an equilibrium price, we have the risk neutral measure $Q$ on this 
filtered space, a probability measure equivalent to the given one. 
And with $B$ any standard $Q$-Brownian motion, the 
exact modelling then is that $S$ is the strong solution of the 
following stochastic differential equation:
$$ dS_t=\varpi\cdot S_t\cdot dt +\sigma \cdot S_t\cdot dB_t, 
\qquad t\in[0,\infty)\, ,$$
or equivalently using It\^o calculus,
$S_t=S_0\exp((\varpi\!-\!\sigma^2/2)t+\sigma B_t)$, for all 
$t$ in $[0,\infty)$.
The positive constant $\sigma$ is the volatility of $S$. The 
specific form of the otherwise arbitrary constant $\varpi$ depends 
on the nature of the security modelled which could be a stock,  a 
currency, a commodity etc. For example,  if $S$ is the price of a 
stock paying  a dividend continuously so as to have constant 
dividend yield $\delta$, then we have $\varpi=r\!-\!\delta$.
%
\grosserAbsatz
{\large\bf 3.\quad Asian options and their equilibrium pricing:}\quad
In the Black-Scholes framework of \S 2, fix 
any time $t_0$ and consider the process $J$  given for 
any time $t$ by:
$$ J(t)=\int_{t_0}^t S_u \, du\, ,$$
The {\it arithmetic-average Asian option\/} written 
at time $t_0$ with maturity $T$ and strike price $K$ is then
the stochastic process on the closed time interval from 
$t_0$ to $T$ paying 
$$ \left( {J(T)\over T-t_0}-K\right)\vbox to 10pt {}^+=\max\left\{
0, {J(T)\over T-t_0}-K\right\}$$
at time $T$ and paying nothing at all other times. As such it is 
a contingent claim on the time interval from $t_0$ to $T$ with 
payoff $ (J(T)/(T\!-\!t_0)-K)^+$.  
\mittlererAbsatz
It is one of the fundamental insights that in the the equilibrium 
framework of the Black-Scholes model, any such contingent claim on 
a risky security has an equilibrium  price which is equal to the 
discounted expectation of its payoff with respect to the risk neutral 
measure conditional on today's information, see 
\cite[\S 22K, (47)]{Duffie88},
\cite[\S 8A]{Duffie92} or [{\bf MR}, Corollary 5.1.1], for instance. 
Applying this {\it arbitrage pricing principle\/}, 
the price $C_t$ of the Asian option at any time $t$ between $t_0$ 
and $T$ is given by the discounted $Q$-expectation conditional on the 
information $\fa_t$ available at time $t$:
$$C_t=e^{-r(T-t)} E^Q\left[\left( 
{J(T)\over T-t_0}-K\right)^+\bigg| \fa_t\right]\, . $$ 
%
However, following [{\bf GY}, \S 3.2], we do not focus on this price, 
but instead we normalize the valuation problem as follows. On factoring 
out the reciprocal of the length $T\!-t_0$ of the time period,
we split the integral $J(T)$ into two integrals, one of which is 
deterministic by time $t$ and the other of which is random. 
We then couple the deterministic integral with the new strike. For 
the random integral, we restart the Brownian motion driving 
the underlying at time $t$, and then using the scaling 
property of Brownian motion, we change time so as to normalize 
its coefficient in the new time scale to two. 
The precise result is the factorization:
$$ C_t={e^{-r(T-t)}\over T-t_0}\cdot {4S_t\over \sigma^2}\cdot
C^{(\nu)}(h,q)\, ,  $$
which reduces the general valuation problem to computing 
$$ C^{(\nu)}(h,q)=E^Q\bigl[(A^{(\nu)}_h-q)^+\bigr], $$
the {\it normalized time-$t$ price\/} of the Asian option.
Herein, $A^{(\nu)}$ is Yor's process
$$
A^{(\nu)}_h
=\int_0^h e^{2(B_w+\nu w)}dw \, , 
$$
and the normalized parameters are as follows:
$$
\nu= {2\varpi\over \sigma^2}-1,\qquad 
h={\sigma ^2\over 4}(T-t),\qquad 
q=kh\!+\! q^*,
$$
where
$$
k={K\over S_t}, \qquad q^*=q^*(t)=
{\sigma^2\over 4S_t}\left( K\cdot(t\!-\!t_0)
-\int_{t_0}^t S_u\, du \right).
$$
\kleinerAbsatz
To interpret these quantities,  
$\nu$ is the {\it normalized adjusted interest rate\/},  $h$ is
the {\it normalized time to maturity\/}, which is 
non-negative, and $q$ is the {\it normalized strike price\/}.   
On a conceptual level, notice that valuing any Asian option in 
this way is reduced to computing a single function $C^{(\nu)}$, 
and that a similar notion of normalized hedging of Asian options 
can be developed along these lines. On a structural level, moreover
notice how $q$ becomes affine linear in the time variable $h$ 
with coefficients $k$ and $q^*$ depending only on  quantities 
known at time $t$.
\goodbreak\grosserAbsatz
%
%
{\large\bf 4.\quad A first reduction of the normalized valuation:}\quad
There is now a dichotomy in computing the normalized time-$t$ price 
$$ 
C^{(\nu)}(h,q)=
E^{Q}\big[ \big( A^{(\nu)}_h-q\big)^{\!+}\big]
$$
of the Asian option according to the normalized strike 
price $q$ being positive~or~not. 
Indeed, if $q$ is non-positive, Asian options lose their option feature. 
Computing  the values  
$C^{\raise-2pt\hbox{$\scriptstyle (\nu)$}}(h,q)$ is straightforward:
\grosserAbsatz
{\bf Lemma:} {\it If $q\le 0$, we have 
$$ C^{(\nu)}(h,q)
=E^{Q}\Big[ A^{(\nu)}_h\Big]-q\, . 
$$ 
with 
$$
E^{Q}\Big[ A^{(\nu)}_h\Big]
= {e\vbox to 7pt{}^{2h(\nu+1)}-1\over 2(\nu\!+\!1)}
\, ,
$$
for any real $\nu$, and it is thus sufficient to compute 
$C^{(\nu)}(h,q)$ if $q>0$.} 
\grosserAbsatz
This can be proved on applying Fubini's theorem, and 
the last expectation is seen to be analytic in $\nu$ 
with its value at $\nu=-1$ equal to $h$. It should be 
mentioned that formulas for all moments of $A^{(\nu)}$
have been derived at various instances over the last 
fifty years, see for example 
\cite[\S 2.4.1, (4.d'), p.33 and Postscript \#3b), p.54]{Yor01}. 
\goodbreak\grosserAbsatz
%
{\large\bf 5.\quad Yor's integral representation for 
Asian option values:}\quad
A closed form for the normalized time-$t$ prices $C^{(\nu)}$ 
of Asian options can be obtained as  a consequence of Yor's 
triple integral representation \cite[(6.e), p.528]{Yor 92}. 
Recall the latter is based on Yor's characterization of the 
law of $A^{(\nu)}$ in \cite[(6.c), p.527]{Yor 92} and so is 
eventually based on the Hartman-Watson theory of \cite{Yor 80}. 
Furnishing a measure for the difficulty of computing normalized 
prices, the precise form of Yor's closed form is as follows
\grosserAbsatz
{\bf Theorem:} {\it For any reals $h$, $q>0$ and $\nu$, we have}  
$$C^{\raise -1.5pt\hbox{$\scriptstyle (\nu)$}} (h,q)
=c_{\nu,h}
\int_0^\infty x^\nu\int_0^\infty 
 e\vbox to 10 pt{}^{-{\scriptstyle (1+x^2)y\over\scriptstyle 2}}
\cdot\Bigl( {1\over y}\!-\! q\Bigr)^{\hskip-2.pt +}
\cdot \psi_{xy}(h)\, dy\, dx\, .$$
\kleinerAbsatz
Herein the function $\psi_a$, for any positive real number $a$, 
is given by the following integral
$$\psi_a(h)
=\int_0^\infty 
 e\vbox to 10 pt{}^{-{\scriptstyle w^2 \over\scriptstyle 2h}}
 e\vbox to 10 pt{}^{-a \cosh(w)}
\sinh(w)\cdot \sin\Bigl( {\pi\over h}w\Bigr) dw\,  , $$
for any $h>0$, and we abbreviate
$$
c_{\nu,h}={1\over \pi \sqrt{2\pi^3 h\, }} 
\, e\vbox to 10 pt{}^{
{\scriptstyle \pi^2\over\scriptstyle 2h}
-{\scriptstyle\nu^2 h\over\scriptstyle 2}}\,  .
$$ 
While Yor's formula seems to require $\nu$ to 
be bigger than at least minus one, it is valid for all real $\nu$.
This is proved in [{\bf SL}, \S 6]. 
\mittlererAbsatz
One purpose of closed form expressions is to provide means for actually 
computing option prices. While Yor's results above are the 
key to many insights into the mathematical structure of $A^{(\nu)}$, 
Yor's formula of the Theorem has a number of structural difficulties
in this regard. First, it involves three integrations with seemingly 
no further structure or further possibilities for simplification. 
Thus, methods for explicitly computing it will necessarily be rather 
complex. However, we have noticed an apparently even bigger obstacle 
to computability. For instance taking $t_0=0$ and $\varpi=r$ equal to $5\%$, 
we compute the factors $c_{\nu,h}$ as follows:
$$ 
\vcenter{
\vbox{\offinterlineskip 
\hrule
\halign{&\vrule \hfil #\vbox to 11pt{}\hfil & \strut \enspace  $\hfil#\hfil$\enspace\cr
&c_{\nu,h} &&
\sigma=20\% &&\sigma=30\%&&\sigma=40\%   &\cr 
\noalign{\hrule}
&\hbox{$T=1$ year}&& 
2\,.\,627\times 10^{213}
&& 2\,.\,265\times 10^{94}&& 4\,.\,816\times 10^{52}&\cr
\noalign{\hrule}
& \hbox{$T=6$ months} && 
7\,.\,686\times 10^{427}
&& 5\,.\,717\times 10^{189}&& 2\,.\,583\times 10^{106}&\cr
\noalign{\hrule}
\noalign{\vskip3pt}
\multispan{9} \hfill \hbox{{\ninebf Table 1.}\enspace 
\ninerm ${\scriptstyle c_{\nu,h}}$ as function of 
${\scriptstyle T}$  and ${\scriptstyle \sigma}$.}\hfill\cr 
} 
}
}
$$
As the examples of \S 8 or \S 24 will show,  normalized prices of 
Asian options are not too big. Yor's formula thus expresses them 
as the product of a big number times a triple integral. The latter 
so has to be very small and must be computed with very high 
accuracies to get reasonably accurate results. The Laplace transform 
approach developed in \cite{Geman+Yor} explicitly for the 
purpose of valuing Asian options was seen to offer a way out 
of these difficulties.
\vskip.8cm
\centerline{\Large\bf II\qquad Laplace transform results}
\vskip.3cm\noindent 
%
{\large\bf 6.\quad The Laplace transform in option valuation:}\quad
Working with continuous functions on the non-negative real 
line of at most exponential growth, the {\it Laplace 
transform\/} $\la(f)$ of any such function $f$ is  
defined by
$$
\la(f)(z)=\int_0^\infty e^{-zx} f(x)\, dx\,,
$$   
for any complex number $z$ in a half-plane contained 
sufficiently deep within the complex right half-plane.
\mittlererAbsatz
The connection of this notion with option valuation in
general and Asian option valuation in particular is as
follows. Fix any option type, like \S 3's European
style Asian call option on a certain stock with price $S$. 
It will depend on a number of parameters, like 
strike price and maturity date.  
At a fixed point in time $t$, consider the family 
which consists of all options on the market with all 
such parameters fixed except maturity dates $M$. In the 
Asian option example thus consider the Asian options 
of all maturities available at time $t$ which have the 
same strike price. In this way regard maturity date $M$ as a 
real variable ranging from $t$ to infinity. The value
of the option thus becomes a function of $M$. 
\mittlererAbsatz
However, it is normalized prices $C^{(\nu)}$ we have to consider
for the Asian option. The normalizations of \S 3 in fact turn 
$C^{(\nu)}$ into a function of, in particular, normalized time 
to maturity.  As a function of maturity date $M$ normalized time 
to maturity is explicitly given by $h(M)=(\sigma^2/4)\cdot(M\!-\!t)$. 
With $M$ from $t$ to infinity, $h(M)$ thus ranges from $0$ to 
infinity. The normalized price of the Asian option so becomes 
a function on the non-negative real line. Call it $f_{AO}$ for 
the sake of emphasis. Recalling from \S 3 how the normalized 
strike price $q$ depends in an affine linear way on normalized 
time to maturity, we then have more precisely
$$
f_{AO}(x)
=E^Q\Big[\big(A_x^{\raise-2pt\hbox{$\scriptstyle (\nu)$}}\!-\! (kx\!+\!q^*)
\big)^+\Big]
$$
for any non-negative real $x$. 
\mittlererAbsatz
At this point we should signal a difficulty with \cite{Geman+Yor}
which we have noticed in Fall 1999. It is not the functions $f_{AO}$ 
this last paper is working with, but the functions $f_{GY,a}$ given 
for any positive real $a$ by
$$
f_{GY,a}(x)
=E^Q\Big[\big(A_x^{\raise-2pt\hbox{$\scriptstyle (\nu)$}}\!-\!a
\big)^+\Big]
$$
for any non-negative real $x$. Notice that the function $f_{AO}$ 
has a non-constant strike price $x\mapsto kx\!+\!q^*$ whereas
any function $f_{GY,a}$ has the constant strike price $a$. So 
the value function of the Asian option $f_{AO}$ is different from 
any function $f_{GY,a}$ which thus is the  value function of a  
{\it non-Asian option\/}. Because of the injectivity of the Laplace 
transform on continuous functions, the Laplace transform 
$\la(f_{AO})$ is different from any Laplace transform 
$\la(f_{GY,a})$ too. The problem at this point so is one of relevance 
of the mathematics for the finance application proper: is there 
a way of relating the valuation of Asian options to the valuation 
of the non-Asian options?   
\goodbreak\grosserAbsatz
%
%
{\large\bf 7.\quad Valuing Asian options using families 
of non-Asian options:}\quad
The basic idea is as follows. Try to reconstruct the 
normalized time-$t$ price 
$$C^{\raise-2pt\hbox{$\scriptstyle (\nu)$}}(h,q)=
E^Q\Big[\big(A_x^{\raise-2pt\hbox{$\scriptstyle (\nu)$}}\!-\! (kx\!+\!q^*)
\big)^+
\Big]\Big|_{x=h}$$
of the Asian option from a family of auxiliary functions whose single
members are unrelated to the problem of valuing the Asian option but 
amenable to the \cite{Geman+Yor} ana\-lysis. Given that the time 
dependency of the normalized strike price poses the problems, simply 
force this strike price to be constant. Thus arrive, for any positive 
real $a$, at \S 6's functions $f_{GY,a}$ on the positive real line  
recalled to be given by
$$ 
f_{GY, a}(x)=E^Q\big[ (A_x^{(\nu)}\!-\!a )^+\big]
$$
for any positive real $x$. These are the functions considered in 
\cite{Geman+Yor}. As we remarked in \S 6  they are the values 
of certain non-Asian options, and taken  individually, they 
cannot be used to value the original Asian option. However, our
finding is that as a whole they allow one to recover the normalized 
time-$t$ price. With the concepts of \S 3 this {\it key reduction\/} 
is made precise in  our following
\mittlererAbsatz
{\bf Lemma:}\quad{\it  If $q=kh\!+\!q^*$ is positive,
computing the normalized time-$t$ price 
$C^{\raise-1.5pt\hbox{$\scriptstyle (\nu)$}}(h,q)$ of the Asian option 
reduces to computing all $f_{GY,a}$ with $a>0$. More precisely, 
$C^{\raise-1.5pt\hbox{$\scriptstyle (\nu)$}}(h,q)$ is obtained by 
choosing the function $f_{GY, kh+q^*}$ and evaluating it at $h$. 
} 
\mittlererAbsatz
A moment's reflection will convince the reader that this is true 
by construction, and meanwhile, the particular case of the Lemma 
we explained to Yor in May 2000 can be found in \cite[pp. 95--96]{Yor01}.  
Notice that our construction works in the more general situation 
where functions of the form $h(y)=f(y,\varphi(y))$ with a known map 
$\varphi$ have to be computed: compute the functions $f$ and 
then intersect with the graph of $\varphi$ to get $h$. And to stress 
the main result again, as one consequence of this technique, we have 
at this point reduced valuing Asian options to valuing the family of 
all non-Asian options.     
\goodbreak\grosserAbsatz
{\large\bf 8.\quad Laplace transforms of the non-Asian 
option values $f_{GY,a}$:}\quad
The key reduction of the preceding \S 7 shows the way to apply 
the \cite{Geman+Yor} Laplace transform to value Asian options. 
Moreover, adopting the notation of \S 3,  consider the family 
of all functions $f_{GY,a}$ with $a>0$ that send any $x>0$ to  
$$ f_{GY,a}(x)=
E^Q\big[(A\raise-1.7pt\hbox{${}_{\raise-1.pt\hbox{$\scriptstyle x$}}^{
\raise-1.4pt\hbox{$\scriptstyle(\nu)$}} $}\!-\!a)^+\big].$$
Recall from \S 7 that its single members are unrelated to valuing the Asian 
option, but that as a whole, the family allows reconstruction of the
normalized time-$t$ price 
$C^{\raise-1.pt\hbox{$\scriptstyle(\nu)$}}(h,q)$ of the Asian option. 
As a first step in actually computing the $f_{GY,a}$, try to compute 
their La\-pla\-ce trans\-form $ F_{GY,a}$ given by 
$$ F_{GY,a}(z)=\int_0^\infty e\vbox to 8pt{}^{ -zx} f_{GY,a}(x)\, dx
=\Laplacetrafo (f_{GY,a})(z).$$
Here the complex number $z$ is to be taken in a half-plane 
sufficiently deep within the right complex half-plane such that
the integrals are finite. The function so obtained is analytic. 
The precise conditions under which these integrals are finite
is part of our description of these 
{\it generalized Geman-Yor Laplace transforms\/} 
in~the~following
\grosserAbsatz
{\bf Theorem:}\quad {\it 
If the normalized strike price $q$ is positive, 
the integrals $F_{GY,a}$ are finite for any complex number 
$z$ with $\re (z)>\max\{0, 2(\nu\!+\!1)\}$, and we have 
$$
F_{GY,a}(z)={D_\nu(a,z)\over z(z\!-\!2(\nu\!+\!1))} $$
where on choosing the principal branch of the logarithm}
$$
D_\nu(a,z)
=
{e\vbox to 6pt{}^{ -{\scriptstyle 1\over\scriptstyle 2a}}
\over a}
\int_0^\infty 
e\vbox to 8pt{}^{-{\scriptstyle x^2\over\scriptstyle 2a}}
x^{\nu+3}
I_{\sqrt{2z\!+\! \nu^2}}\Big( {x\over a}\Big)dx\, .$$ 
\kleinerAbsatz
Herein $I_\mu$ is the modified Bessel function with complex order 
$\mu$, as discussed in \cite[Chapter 5]{Lebedev}. 
Generalizing \cite[(3.9), p.363]{Geman+Yor}, the integral $D_\nu(a,z)$ 
can be expressed using the confluent hypergeometric function $\Phi$ 
discussed in  \cite[Chapter 9]{Lebedev} as follows
\grosserAbsatz
{\bf Corollary:} \quad{\it  For any complex $z$ with 
$\re(z)>\max\{0,2(\nu\!+\!1)\}$, we have
$$
D_\nu(a,z)=
{\Gamma\Big( {\displaystyle \nu\!+\!4\!+\!\mu(z)\over\displaystyle 2}\Big)
\over \Gamma\big(\mu(z)\!+\!1\big)}\cdot
  \Phi\bigg( {\nu\!+\!4\!+\!\mu(z)\over 2 }, 
              \mu(z)\!+\!1; {1\over 2a}\bigg)
\cdot {         
(2a)\vbox to 9pt{}^{{\scriptstyle \nu+2-\mu(z)\over\scriptstyle 2}}
\over
e\vbox to 6pt{}^{ {\scriptstyle 1\over\scriptstyle 2a}}
}\, 
$$  
\vskip3pt\noindent
on setting $\mu(z)=\sqrt{2z\!+\!\nu^2\, }$.}
\grosserAbsatz
We again stress that these Laplace transforms
are not those of the Asian option price, but rather the Laplace
transforms of the prices of auxiliary options. To obtain Asian
option prices, we have to invert these Laplace transforms and 
then proceed using \S 7~Lemma; formally speaking:
$$
C^{(\nu)}(h,q)=\la^{-1}\left( F_{GY,q}\right)(h)\, . 
$$  
In full mathematical generality,  analytic inversion has been 
achieved in \cite{Asia}, and we come back to this in \S 23. 
Numerical inversions have also been accomplished. For example  
\cite{Fu+Madan+Wang} computed the following seven cases as 
reproduced in \cite[Table 7.1]{Dufresne00}
$$
\vcenter{
\vbox{\offinterlineskip 
\hrule
\halign{&\vrule #& \strut \enspace  $\hfil#\hfil$\enspace\cr
&\hbox{\rm Case}&& r&&\sigma && T && S_0 && \nu && 2C^{(\nu)}&\cr
\noalign{\hrule}
&1&& 2\% && 10\% && 1 && 2.0&&  3&& 0.056&\cr 
\noalign{\hrule}
&2&& 18\% && 30\% && 1 && 2.0&& 3&& 0.219&\cr
\noalign{\hrule}
&3&& 1.25\% && 25\% && 2 && 2.0&&  -\,0\,.\,6&& 0.172&\cr
\noalign{\hrule}
&4&& 5\% && 50\% && 1 && 1.9&& -\,0\,.\,6&& 0.194&\cr
\noalign{\hrule}
&5&& 5\% && 50\% && 1 && 2.0&&-\,0\,.\,6&& 0.247&\cr
\noalign{\hrule}
&6&& 5\% && 50\% && 1 && 2.1&&-\,0\,.\,6&& 0.307&\cr
\noalign{\hrule}
&7&& 5\% && 50\% && 2 && 2.0&&-\,0\,.\,6&& 0.352&\cr
\noalign{\hrule}
\noalign{\vskip5pt}
\multispan{15}  \hfill 
\hbox{\ninerm{\ninebf Table 2.}\enspace 
\ninerm Prices ${\scriptstyle 2C^{(\nu)}}$ for ${\scriptstyle K=2.0}$ 
and}\hfill\cr \noalign{\vskip2pt}\multispan{15} 
\hfill\hbox{\ninerm ${\scriptstyle t_0=t=0}$ using numerical Laplace 
inversion.}\hfill\cr 
} 
}}
$$
Obtaining our two results proceeds in two steps. In a first 
probabilistic step, the arguments of \cite{Geman+Yor} are adapted 
to compute the modified Geman-Yor transforms $F_{GY,a}$. We
give in Part~III a new exposition of the argument incorporating 
the tutorials and kind suggestions of Yor. The key idea is to 
factorize the geometric
Brownian motion of the underlying over a Bessel process of index
$\nu$. Pertinent notions are discussed in \S 9. This makes time
stochastic in such a way that Yor's process $A^{(\nu)}$ now
takes the double role of both a stochastic clock and a control
variable for the Asian option. At first sight, this appears to
complicate the original valuation problem. However, this double
role of $A^{(\nu)}$ is especially suited to the
Laplace transform. Indeed,  in contrast to the 
situation for the Asian option, the strike price $a$ of the 
non-Asian option with value function $f_{GY,a}$ is independent 
of time. This makes it possible to reduce the computation of the 
Laplace transform $F_{GY,a}$ to the following  problem: obtain 
explicit expressions for the Bessel semigroup of index $\nu$ and 
for a certain conditional expectation involving first passage 
times of Yor's process $A^{(\nu)}$.
\mittlererAbsatz
However, such results are available for both concepts only if the 
index $\nu$ is non-negative, an assumption which is explicit in 
\cite[\S2]{Geman+Yor}. This non-negativity condition, however, 
translates into the condition that the risk-neutral drift is not 
less than half the squared volatility. Unfortunately, this places 
restrictions on the financial applicability of the result. For 
example, if volatility is 30\%, the arguments of \cite[\S 3]{Geman+Yor}
are not valid if the difference between the riskfree rate and 
the dividend yield is less than 4.5\%, and then we do not have the 
Laplace transforms of the non-Asian options either. Worse yet, 
the greater is the volatility, the greater is this range of parameters 
in which we do not have the Laplace transforms. Unfortunately, it is
precisely due to high volatility that Asian options are used in the 
first place.
\mittlererAbsatz
We have been a bit disconcerted by these findings. Luckily, however,
we found that in particular those results of Table~2 where $\nu$ is 
negative were reproduced in \cite{Dufresne00} using an alternative 
approach. And corroborating Andr\'e Weil's \cite[p.457]{Weil} dictum 
that `theorems are proved by those who believe in them', we are now 
able to discuss in the sequel three different ways for stablishing 
the Theorem and its Corollary for arbitrary real risk neutral 
drifts~$\nu$.
\mittlererAbsatz
The first of these, as discussed in Part IV, is inspired by \cite{Yor 92},
and we think Yor could have given this argument had he been aware of
the financial motivation for letting the parameter $\nu$ be negative.
In fact, the key idea 
is to try to bypass the difficulties of Bessel processes of negative 
index $\nu$ by Girsanov transforming to the simpler Bessel processes 
of index zero.This theme is developed also in the third approach 
discussed in Part VI. Here, the idea is not to 
Girsanov transform Bessel processes derived from the geometric Brownian
motion driving the underlying. Instead, Girsanov transform this geometric
Brownian motion itself; an idea we distilled from \cite{Yor 80}. 
The effect of this is that zero drift Bessel processes enter right from 
the beginnning, and the result is a uniform argument based on these 
most natural Bessel processes. 
\mittlererAbsatz     
In comparison to these two approaches, our second approach discussed 
in Part V seems somewhat novel. The idea is to combine stochastic 
methods with complex analytic ones.  The net effect here is that 
with an input of some standard result from the latter area like 
the identity theorem, it is not Bessel processes which now have to 
be dealt with but Brownian motion. For this it is  moreover not 
required to work on the process level, as it is sufficient to work 
on the level of expectations. So this second approach seems to be 
an example of a rather promising methodology for solving problems, 
which is to systematically enhance stochastic techniques 
with complex analytic ones.
\vskip1cm
\centerline{\Large\bf III\quad Laplace transforms if $\nu\ge 0$}
\vskip.3cm\noindent
%
%
{\large\bf 9.\quad Preliminaries on Bessel processes:}\quad
As a preliminary to establishing  the Laplace transforms of \S 8, 
this section collects a number of pertinent facts from the theory 
of Bessel processes. This theory is patterned after the example of 
the Bessel processes of integer dimension $\delta\ge 2$, which are  
defined by taking the Euclidean distance from the origin of a Brownian 
motion in dimension $\delta$. Applying It\^o's Lemma, their 
infinitesimal generator 
${\text{\script{A\hspace*{4.85pt}}}}$ is seen to be given by
$${\text{\script{A\hspace*{4.85pt}}}}f(x)=
{1\over 2}f''(x)+{2\nu\!+\!1\over 2x}f'(x)\, , 
$$
for any function $f$ in $C^2_b({\bf R}_{>0})$. This notion makes 
sense for any real number $\delta$, and the real-valued diffusion 
associated to ${\text{\script{A\hspace*{4.85pt}}}}$ using the 
Volkonskii construction, see for instance \cite[Theorem 4.3.3, p.91]{Knight},
is the {\it Bessel processes\/} 
$R^{ \raise-1.5pt\hbox{$\scriptstyle (\nu)$}}$ on $[0,\infty)$ 
with {\it index\/} $\nu=(\delta/2)-1$. 
While Bessel processes of non-negative indices $\nu$ stay positive if
started with a positive value at time zero, Bessel processes of 
negative indices $\nu$ develop some pathologies. As explained in
\cite[XI \S 1]{Revuz+Yor}, in this case they hit zero. If $-1<\nu<0$, 
they are thereupon instantantanously reflected and never become negative. 
For $\nu=-1$ they continue at zero.  
\mittlererAbsatz
This matters for the second way of defining Bessel processes
of arbitrary dimension $\delta$. 
In fact, the focus here is on squares of Bessel processes. Applying 
It\^o's Lemma, they are the continuous strong solutions of the 
stochastic differential equation
$$ d\rho_t=2\delta dt +2\sqrt{|\rho_t|\, }\, dB_t\, , \qquad \rho_0=1$$
\cite[XI, \S 1]{Revuz+Yor}. These stochastic differential equations 
make sense for any real number $\delta$ and have a unique continuous 
strong solution also if $\delta$ is smaller than two. The so obtained
processes are studied in \cite[\S3]{Goeing-Yor}. For non-negative 
indices $\nu$, they coincide with the squares of the corresponding 
Bessel processes of index $\nu$. They develop some pathologies for 
negative indices $\nu$. In this case, they hit zero if started with 
a positive value at time zero, and if  $\nu<-1$, they even continue 
negative. Notice that in such situations their square roots are purely
imaginary, and so cannot coincide with any of the Bessel processes
constructed above. However, these two ways of extending the notion of Bessel
processes do coincide for processes started at time zero at a 
positive value up to the first time zero is hit. This  
essentially is the reason behind the following {\it Lamperti identity\/},
which may nevertheless be surprising
\grosserAbsatz
{\bf Lemma:}\quad{\it For the index $\nu$ any real number, we have
the factorization:
$$
e\vbox to 8pt{}^{B_t+\nu t}
=R^{(\nu)}\big( A_t^{(\nu)}\big).$$
for any $t>0$, where 
$ A_t^{(\nu)}=\int\nolimits_0^t e\vbox to 8pt{}^{2(B_w+\nu w)}dw$
is Yor's process.}
\mittlererAbsatz
For $\nu\ge0$ a proof is given in \cite[\S 2]{ Yor 92a} while the 
general case is now contained as exercise XI~(1.28), p.452 in the 
third edition of \cite{Revuz+Yor}. We are indebted to Yor for this 
and for kindly supplying us with the following argument. 
\mittlererAbsatz
The idea is to apply the It\^o rule to the square $Z_t$ of 
$Y_t=\exp(\nu t\!+\!B_t)$ to get 
$$  
Z_t=2(\nu\!+\!1)\int_0^t Z_w dw + 2\int_0^t Z_w dB_w \, .$$
Time change the process using the inverse function
$\tau(t)=\inf\{ u| \hbox{$\int_0^u Z_w \, dw >t$} \}$
to Yor's process $A^{(\nu)}$ to get
$$
Z_{\tau(t)}= 2(\nu\!+\!1)t+2\int_0^{\tau(t)} Z_w dB_w \, .$$
To interpret the stochastic integral in this sum, apply the basic 
time change formalism for stochastic processes as in 
\cite[\S 8.5]{Oeksendal} to obtain  
$$\int_0^{\tau(t)} Z_w dB_w
  =\int_0^t  Z_w \sqrt{ {\tau'(w)}\, }dW_w\,  
$$
where $W_t$ is defined as the stochastic integral
$   
W_t=\int_0^{\tau(t)}\sqrt{ Z_w\, }\, dB_w  
$ 
and is a Brownian motion. Using the inverse function theorem of 
calculus, the derivative of $\tau$ is equal to the reciprocal of 
the derivative with respect to time of Yor's process $A^{(\nu)}$ 
at time $w$. Hence it is equal to the reciprocal of $Z_w$. On 
substitution we so identify the time changed process $Z$ as a 
continuous solution to the stochastic differential 
equation for the square of the Bessel process of index $\nu$:
$$ 
Z_{\tau(t)}= 2(\nu\!+\!1)t+2\int_0^t \sqrt{Z_w\, }\, dW_w\, . 
$$
Using the uniqueness of the solution of these stochastic differential 
equations, the time-changed process $Z$ is the square of the 
Bessel process of index $\nu$. Reversing the time change, this 
translates into
$$
Y_t^2=(R_t^{(\nu)})^2( A_t^{(\nu)}). $$
To establish the identity of the Lemma, we have to take square roots. 
This is not a problem if $\nu$ is non-negative since then the  
Bessel process takes non-negative values only. It does pose a problem 
if $\nu$ is negative. In this case, however, recall that the 
Bessel process starts at time zero with a positive value. Since it is 
continuous by hypothesis, it will stay positive until it first hits
zero at time $t^*>0$. Since the process $A^{(\nu)}$ starts at
zero at time zero, there is a latest point in time $t^{**}$, infinity 
admitted, such that $A^{\raise-1.5pt\hbox{$\scriptstyle (\nu)$}}$ 
is smaller than $t^*$ at all points in time $t$ smaller than $t^{**}$. 
Thus we have the required identity at least for all points in time $t$ 
smaller than $t^{**}$. Now $Y_t$ is never zero. Since the processes 
on both sides of the identity are continuous in time, $t^{**}$ must be 
infinity, and the proof is complete. 
%
%
\goodbreak\grosserAbsatz
{\large\bf 10.\quad Computing Laplace transforms if $\nu\ge 0$:}\quad
This section is the first step in the proof of the integral 
representation of \S 8~Theorem for the Laplace transform 
$$ F_{GY,a}(z)=\int_0^\infty e\vbox to 8pt{}^{ -zx} f_{GY,a}(x)\, dx
=\Laplacetrafo (f_{GY,a})(z), $$
where, with the concepts of \S 3 and \S 6, we have
$ f_{GY,a}(x)=
E^Q[(A\raise-1.7pt\hbox{${}_{\raise-1.pt\hbox{$\scriptstyle x$}}^{
\raise-1.4pt\hbox{$\scriptstyle(\nu)$}} $}\!-\!a)^+] 
$, 
for any positive real numbers $a$ and $x$. We now explain why
one needs to restrict the probabilistic arguments of 
\cite{Geman+Yor} and apply them {\it mutatis mutandis\/} in 
order to arrive at the following
\grosserAbsatz
{\bf Lemma:}\quad {\it The assertions of \S 8~Theorem are valid
if $\nu=2\sigma^{-2}\varpi\!-\!1\ge 0$.}
\mittlererAbsatz 
We are very indebted to Yor for correspondence and discussions 
about this result, and are very grateful for his kind support. 
In the sequel, we want to indicate the key steps of the proof 
following \cite{Geman+Yor}, while trying to incorporate his 
suggestions. Hopefully, no pitfalls have remained undetected. 
\goodbreak\mittlererAbsatz 
The basic idea is to make time stochastic using 
the Lamperti identity 
$$
e\vbox to 8pt{}^{ \nu w+B_w}
= R^{(\nu)}\big( A_w^{(\nu)}\big) 
$$
for all positive real numbers $w$ as has been discussed in the 
preceding section. Here, $R^{(\nu)}$ is the Bessel 
process of index $\nu$ with 
$R^{\raise-1pt\hbox{$\scriptstyle (\nu)$}}(0)=1$. 
On applying this Lamperti identity, 
$A^{\raise-2pt\hbox{$\scriptstyle (\nu)$}}$
has the double role 
of both control variable and stochastic~clock. That the ``strike price'' $a$ 
is independent of time now becomes essential. It makes possible to 
transcribe the condition on the control variable  to be bigger than $a$ 
as the inverse time change  
$$\tau_{\nu,a}=\inf\{ u\, | \, A^{\raise-2pt\hbox{$\scriptstyle (\nu)$}}_u>a\}
$$
for the stochastic clock. This is the key idea for obtaining
the representation  
$$
f_{GY,a}(w)
   =E^Q\bigg[ 
     {e\vbox to 8pt{}^{ 2(\nu+1)[w-\tau_{\nu,a}]^+}-1
       \over 2(\nu\!+\!1)}\cdot (R^{(\nu)}_a)^2
       \bigg],
$$
for all $w>0$. Indeed, fix any positive real 
number $x$,  and consider the process $A^{(\nu)}$ at $x$ on the 
set of all events where $\tau_{\nu,a}$ takes values 
less than or equal to $x$. Break the integral defining 
$A^{\raise-2pt\hbox{$\scriptstyle (\nu)$}}(x)$  
at  $\tau_{\nu,a}$.  The first summand then is 
$A^{\raise-2pt\hbox{$\scriptstyle (\nu)$}}$ at 
time $\tau_{\nu,a}$ and so is equal to $a$. In the second summand,
restart the Brownian motion in the exponent of the 
integrand at $\tau_{\nu,a}$ shifting the variable of integration 
accordingly. The second integral then is the product of 
$\exp(2\cdot(B(\tau_{\nu,a}) \!+\!\nu\tau_{\nu,a}))$
times 
$A^{\raise-2pt\hbox{$\scriptstyle (\nu)$}}$ at $x\!-\!\tau_{\nu,a}$,
by abuse of language after having applied the Strong Markov property. This 
last process is such that it is independent of the information at 
time $\tau_{\nu,a}$. Unravelling the definition of $\tau_{\nu,a}$, 
the first above factor so is the square of the Bessel process 
$R^{\raise-2pt\hbox{$\scriptstyle (\nu)$}}$ 
at time $a$.  Now taking the expectation conditional on the 
information at $\tau_{\nu,a}$, we thus get:
$$ E^Q\Bigl[  \bigl( A^{(\nu)}_x-a\bigr)^+\, \Big|\,
 \Fgeschwungen _{\tau_{\nu,a}}\Bigr]
=\bigl( R^{(\nu)}_a\bigr)^2\cdot  
E^Q\Bigl[A^{(\nu)} _{[x-\tau_{\nu,a}]^+}
\Big].$$
On substitution for the expectation of $A^{(\nu)}(w)$
from \S 4~Lemma  or using  \cite[\S 4]{Yor 92}, 
the required expression for $f_{GY,a}$ follows.
\mittlererAbsatz 
At first sight this appears to complicate the problem.
However, it is just what is especially suited to 
the Laplace transform $F_{GY,a}$  of $f_{GY,a}$ now given by:
$$
F_{GY,a}(z)=\int_0^\infty e\vbox to 8pt{}^{-zw}
E^Q\bigg[ 
{e\vbox to 8pt{}^{ 2(\nu+1)[w-\tau_{\nu,a}]^+}-1
\over 2(\nu\!+\!1)}\cdot (R^{(\nu)}_a)^2
\bigg]
\, dw\, . 
$$
Still, for computing this integral one wants to interchange   
the Laplace integral with the expectation $E^Q$. If $z$ is real, it 
seems best to follow Yor's proposal for justifying this. 
Indeed, with the integrand of the double integral in question 
positive and measurable, apply Tonelli's theorem now but justify 
only in a later step that any of the resulting integrals are finite. 
The case of a general argument $z$ is reduced to this case 
by considering 
the absolute value of the integrand, and the result 
is the identity 
$$ F_{GY,a}(z)=
{1\over z(z\!-\!2(\nu\!+\!1))}
E^Q\Big[ e\vbox to 8pt{}^{-z\tau_{\nu,a}} (R^{(\nu)})^2 \Big]
$$
of measurable functions for any complex number $z$ 
with $\re(z)> 2(\nu\!+\!1)$.
The idea for identifying the expectation in the numerator 
as $D_\nu(a,z)$ then is to condition on the Bessel process 
to obtain 
$$
D_\nu(a,z)=\int_0^\infty \hskip-4ptx^2\,
E^Q\Big[ e\vbox to 8pt{}^{-z\tau_{\nu,a}}\Big| R^{(\nu)}_a=x \Big]
\, p_{\nu,a}(1,x)\, 
dx\, , $$
where $p_{\nu,a}$ is the  time-$a$ semigroup density of the Bessel
process of index $\nu$ starting at $1$ at time zero. Following 
\cite[p.362]{Geman+Yor}
we make this integral explicit by making the single factors of 
its integrand explicit. For this, work with the results recalled
in \hbox{\cite[\S2]{Geman+Yor}}. With respect to the conditional 
expectation factor, under the hypothesis $\nu\ge 0$, Yor has 
computed it at positive real arguments $z$ in 
\cite[Th\'eor\`eme 4.7, p.80]{ Yor 80} 
(see also \cite[Lemma~2.1 and Proposi\-tion~2.6]{Geman+Yor}).
Using analytic continuation, the validity of his result can be 
seen to extend to the arguments $z$ in the right half-plane required
in the present situation. This then gives for the conditional 
expectation factor in $D_\nu(a,z)$ the following expression 
as a quotient of $I$-Bessel functions:
$$ 
E^Q\left[ e^{-z\cdot \tau_{\nu,a}
}\Big\vert R^{(\nu)}(a)=w\right]
={I_{\sqrt{2z+\nu^2}}\over I_\nu}\left( {w\over a}\right).
$$ 
Explicit expressions for the Bessel semigroups $p_{\nu,a}$ 
have been known for $\nu>-1$ for some time, see \cite[(4.3), p.78]{Yor 80} 
or \cite[Proposition 2.2]{Geman+Yor}.    
The density $p_{\nu,a}(1,w)$ of the time-$a$ Bessel semigroup with 
index $\nu$ and  starting point $1$ is 
$$ p_{\nu, a}(1,w)=
{w^{\nu+1}\over a}\cdot e^{-{\scriptstyle  1+w^2\over\scriptstyle 2a}}
\cdot I_\nu\left({w\over a}\right)
.$$
Nothing seemed to have been known about such densities 
if $\nu\le -1$ before \cite{Goeing-Yor}. However, the results 
proved there for $\nu<-1$ still need to be handled with care 
as it will be explained in \S 14. The upshot is that, in 
accordance with \cite[\S2]{Geman+Yor}, the above 
decomposition of $D_\nu(a,z)$ seems to give explicit results 
without further qualifications only if if $\nu\ge 0$. 
Then, however, we have the required result
$$
D_\nu(a,z)
=
{e\vbox to 8pt{}^{ -{\scriptstyle 1\over\scriptstyle 2a}}\over a}
\int_0^\infty 
e\vbox to 8pt{}^{-{\scriptstyle x^2\over\scriptstyle 2a}}
x^{\nu+3}
I_{\sqrt{2z\!+\! \nu^2}}\Big( {x\over a}\Big)dx\, . 
$$
There is a further technical point to be taken care of herein: choose 
the principal branch of the logarithm to define the square root on 
the complex plane with the non-positive real line deleted.
\mittlererAbsatz
To complete the Tonelli argument proposed 
to us by Yor and to complete the proof, we have to establish the 
finiteness of this integral for any fixed complex number $z$  
with $\re (z)>2(\nu\!+\!1)$. 
This is a consequence of \S 11~Proposition's convergence analysis 
of these integrals, and granting this result, the proof of the
Lemma is complete.
\goodbreak\grosserAbsatz
{\large\bf 11.\quad Integrability analysis:}\quad
In terms of the concepts of \S 3, this section studies finiteness 
of the integrals
$$
D_\nu(a,z)
=
{e\vbox to 8pt{}^{ -{\scriptstyle 1\over\scriptstyle 2a}}\over a}
\int_0^\infty 
e\vbox to 8pt{}^{-{\scriptstyle x^2\over\scriptstyle 2a}}
x^{\nu+3}
I_{\sqrt{2z\!+\! \nu^2}}\Big( {x\over a}\Big)dx\,  
$$  
for any real $a>0$ in terms of their complex parameters $\nu$ and $z$. 
The precise result is the following
\grosserAbsatz
{\bf Proposition:} {\it Let $\varepsilon\ge0$ be any real number. 
If $|\im(\nu)|\le \varepsilon$, the integrals $D_\nu(a,z)$ are finite 
for any complex $z$ with real part $\re(z)>2\varepsilon^2$.}
\grosserAbsatz
The Proposition depends on the following result about 
the complex square root associated to the principal branch of the 
complex logarithm.
\grosserAbsatz
{\bf Lemma:} {\it Let $\varepsilon\ge 0$ be any real number. For 
any complex $\nu$ with $|\im(\nu)|\le \varepsilon$ we then have 
$$
\re(\sqrt{2z+\nu^2\, }\, )>|\re(\nu)|
$$
for any complex $z$ with $\re(z)>2\varepsilon^2$.}
\goodbreak\grosserAbsatz
A proof of the Lemma based on a close analysis of the 
square root can be found in \cite[\S10]{carr-sch}. 
To prove the Proposition, finiteness of $D_\nu(a,z)$ 
under the conditions of the Proposition follows by 
combining the above square root lemma with the 
asymptotic behaviour of the Bessel function factor 
of its integrand near the origin and towards infinity.
Indeed, setting 
$$
\mu=\sqrt{2z\!+\!\nu^2\, }\, ,
$$
from \cite[\S 5.11]{Lebedev} recall that $I_\mu$ is a continuous 
function on the positive real line in particular whose asymptotic 
behaviour for large real arguments is 
$$
I_\mu(\xi)\approx {e^{\scriptstyle \xi}\over \sqrt{2\pi \xi\, }}
\qquad 
\hbox{for $\xi\rightarrow \infty$}\, . 
$$
Hence the factor $\exp(-x^2/(2a))$ dominates the asymptotic 
behaviour of the integrand of $D_\nu(a,z)$ with $x$ to infinity,
whence its integrability away from the origin. On the other hand, 
from \cite[\S 5.7]{Lebedev} we have for real arguments near zero 
$$
I_\mu(\xi)\approx {\xi^\mu\over 2^\mu\Gamma(1\!+\!\mu)}\qquad 
\hbox{for $\xi\downarrow 0$}
\, .$$
Thus, if the real part of $\mu\!+\!\nu\!+\!4$ is positive, or 
equivalently, if we have 
$$ 
\re(\mu)>-(\re(\nu)\!+\!4),
$$
no integrability problems arise for $x$ near the origin. Under the 
conditions of the Proposition, on the other hand, the above Lemma
gives
$$ 
\re(\mu)>|\re(\nu)|.$$
Since this last inequality implies the former, the proof of the 
Proposition is complete. 
\goodbreak
\vskip.9cm
\centerline{\Large\bf IV\quad Laplace transforms in the general case:} 
\vskip.06cm
\centerline{\Large\bf using Girsanov transforms of Bessel processes}
\vskip.3cm\noindent
{\large\bf 12.\quad Further preliminaries on Bessel processes:}\quad
Our first way of extending the results of \S 10  to negative indices 
$\nu$ requires further preliminaries on Bessel processes 
from \cite[\S 2]{Yor 92}. Recall that Bessel processes where 
the index $\nu$ is any real number are the real-valued diffusions 
whose infinitesimal generators 
${\text{\script{A\hspace*{4.85pt}}}}$ are given by
$${\text{\script{A\hspace*{4.85pt}}}}f(x)=
{1\over 2}f''(x)+{2\nu\!+\!1\over 2x}f'(x)\, , 
$$
for any function $f$ in $C^2_b({\bf R}_{>0})$. 
To describe the law $P_{\mu,u}$ on $C({\bf R}_{\ge 0},{\bf R}_{\ge 0})$ of  
$R^{(\nu)}$ if this process starts at the non-negative real $u$,  
let $\rho$ be the canonical process on $C({\bf R}_{\ge 0},{\bf R}_{\ge 0})$;
recall it operates as evaluation map: $\rho_a(f)=f(a)$. If
${\text{\script{R\hspace*{4.85pt}}}}$ is the canonical filtration
with ${\text{\script{R\hspace*{4.85pt}}}}_a$ equal to the sigma algebra
generated by the $\rho_s$ with $s\le a$, we then have the mutual
absolute continuity relation
\grosserAbsatz
{\bf Lemma:}\quad {\it If the Bessel process of any real index
$\nu$ is started at any non-negative real $u$, its law is 
realated to that of the zero index Bessel process started at $u$
as follows 
$$
{P_{\mu,u}}_{ \big|\,  {\scriptR_{a}} \cap\{ a<T_0\}}
= 
\bigg({\rho_a\over u}\bigg)^\mu 
\exp\bigg(\!\! -{\mu^2\over 2}
\int_0^a {ds\over \rho_s^2}\bigg) 
{P_{0,u}}_{ \big|\,  \scriptR_a}\, , 
$$
where $T_0$ is the first passage time of $\rho$ to zero.}
\grosserAbsatz
This is proved as an application of Girsanov's theorem by exchanging 
drifts in the stochastic differential equation of \S 9, and has
the following
\grosserAbsatz
{\bf Corollary:}\quad {\it For any complex $z$ with positive real part,
and any non-negative real $r$, 
$$
E^{0}_u\bigg[ 
\exp\bigg(\!\! -z 
\int_0^a {ds\over \rho_s^2}\bigg) 
\,  \Big|\, 
\rho_a=r\bigg]=
{I_{\sqrt{2z\,}}\over I_0}\bigg({ur\over a}\Bigg)
\, ,
$$
where the expectation is taken with respect to the law $P_{0,u}$.}
\grosserAbsatz
The Corollary is proved in two steps. If $z$ is any non-negative real,   
$T_0=\infty$, and we have the explicit expressions for the  densities 
$p_{\mu,a}(u,r)$ of the of the Bessel semigroups  already encountered 
in \S 10 and related to the law via $P_{\mu,u}(a,dr)=p_{\mu,a}(u,r)\, dr$:
$$
p_{\mu,a}(u,r)=
\Big({r\over u}\Big)^\mu\, {r\over a}\, 
\exp\bigg(\!\! -{1\over 2a}(u^2\!+\!r^2)\bigg)
I_\mu\bigg( {ur\over a}\bigg), 
$$
for any non-negative reals $u$, $a>0$ and $r$. And so the Corollary
follows on taking expectations in the absolute continuity 
relation of the Lemma. Observing that both sides of the identity
to be proved are analytic functions in $z$ on the right half plane,
the general case then follows by analytic continuation as a second 
step.
\goodbreak\grosserAbsatz
{\large\bf 13.\quad First proof of the Laplace transform using 
Bessel processes:}\quad
Resuming the discussion of \S 10,  we are now able to complete the 
proof of \S 8~Theorem in the way it might have been envisaged by Yor: 
based on a careful analysis of Bessel processes. Recall that we still 
need explicitly compute for negative normalized risk-neutral drifts 
$\nu$, the risk-neutral expectations 
$$
E^Q\Big[ e\vbox to 8pt{}^{-z\tau_{\nu,a}}\big( R_a^{(\nu)}\big)^2\Big] 
$$
where $a>0$ and $\re (z)$ is positive and sufficiently big in 
particular. From the time change part of the argument in \S 9, 
recall that $\tau_{\nu,a}$ as the inverse time change of the process 
$A^{\raise-2pt\hbox{$\scriptstyle (\nu)$}} $ 
at time $a$ is given by
$$
\tau_{\nu,a}=\int_0^a{ds\over\big(R_s^{(\nu)}\big)^2}
\, . $$
Again using \S 9's Lamperti relation 
$$
e\vbox to 8pt{}^{B_w+\nu w}=R^{(\nu)}\big( A_w^{(\nu)}\big)\, , 
$$
which is valid for $w\ge 0$, the point now is that $a$ is smaller 
than the first passage time to zero $T_0$ of the Bessel process  
$R^{\raise-2pt\hbox{$\scriptstyle (\nu)$}} $. Put differently,
this Bessel process which starts in $1$ at time zero still 
is positive at time $a$. Using the absolute continuity relation of 
\S 12~Lemma,  
$$
E^Q\bigg[\big( R_a^{(\nu)}\big)^2
 e\vbox to 8pt{}^{-z\tau_{\nu,a}}\bigg] 
=
E^0_1\bigg[
\rho\vbox to 7pt{}^2_a 
\,e\vbox to 8pt{}^{-z {\scriptstyle\int\limits_0^a}
{\scriptstyle ds\over\scriptstyle\rho_s^2}}
\cdot
\rho\vbox to 7pt{}^\nu_a 
\,e\vbox to 8pt{}^{-{\nu^2\over 2}{\scriptstyle\int\limits_0^a}{
\scriptstyle ds\over\scriptstyle \rho_s^2}}
\bigg]
\, .
$$  
Conditioning on the index-$0$ Bessel process thus gives 
$$
E^Q\bigg[ e\vbox to 8pt{}^{-z\tau_{\nu,a}}\big( R_a^{(\nu)}\big)^2\bigg] 
=
\int_0^\infty 
E^0_1\bigg[
e\vbox to 8pt{}^{-\big(z+{\nu^2\over 2}\big) 
\int\limits_0^a{ds\over\rho_s^2}}\, \bigg| \, \rho_a=\rho\bigg]\,
\rho^{\nu+2}\, p_{0,a}(1,\rho)\, d\rho
\, , $$
With the explicit form of the semigroup recalled in \S 12,     
applying \S 12~Corollary gives
$$
E^Q\bigg[ e\vbox to 8pt{}^{-z\tau_{\nu,a}}\big( R_a^{(\nu)}\big)^2\bigg]
=
{1\over a}
e\vbox to 9pt{}^{-{\scriptstyle 1\over\scriptstyle 2a}}
\int_0^\infty I_{\sqrt{ 2z+\nu^2\, }}\bigg( {\rho\over a}\bigg) 
\, \rho\vbox to 8pt{}^{\nu+3}
\, e\vbox to 9pt{}^{-{\scriptstyle \rho^2\over\scriptstyle 2a}}
\, d\rho
$$ 
as desired. 
An application of \S 11~Proposition then shows finiteness of 
this integral if $z_0=\re(z)$ is positive and bigger than 
$2(\nu\!+\!1)$, and the first proof of \S 8~Theorem is complete. 
Note that meanwhile the argument of\cite[pp. 97--99]{Yor01} 
seems to corroborate our statement that this proof is very much in the 
spirit of Yor.    
\goodbreak\grosserAbsatz
{\large\bf 14.\quad Vista on the use of Bessel processes:}\quad
The strategy of the preceding argument is to bypass the difficulties 
brought about by Bessel processes with negative indices by Girsanov 
transforming to Bessel processes of index zero. One could
ask about a direct attack in the spirit of \S 10. This would require
possession of explicit expressions for both the densities of the respective 
Bessel semigroups and the pertinent conditional expectations of
$\exp(-z\tau_{\nu,a})$. The recent work of G\"oing-Jaeschke 
and Yor in \cite{Goeing-Yor} now in particular provides certain 
analytic expressions for the densities. However, it is not Bessel 
processes as considered in \S 12 which are studied there, but rather
processes obtained as strong solutions to the stochastic differential 
equations of squared Bessel processes, as mentioned in \S 9. 
Recall that for negative indices, the latter processes become 
negative and their square roots, which should give the Bessel 
processes, then are purely imaginary. Still, 
the two notions of Bessel processes thus obtained coincide 
on their respective positive range. There, we have Bessel 
semigroup densities based on those derived in 
\cite[\S 3 Proposition 2, p.21]{Goeing-Yor}. For indices
$\nu<-1$ and time-$0$ starting values $y>0$, they are given by
$$
p_{\nu,t}(x,y)
=
h(x,y,\delta,t)
\, e\vbox to 10pt{}^{\scriptstyle y-x\over \scriptstyle 2t}
\int_0^1
{(1\!-\!w)^{2(\mu-1)}\over w^\mu}
\, e\vbox to 12pt{}^{ 
{\scriptstyle 1\over\scriptstyle 2t}
\big(x w -{\scriptstyle y\over \scriptstyle w}\big)}
\, dw
$$
defining
$$
h(x,y,\delta,t)={1\over \Gamma^2(\mu\!-\!1)}
\, {(xy)\vbox to 7pt{}^{ \mu-1 } 
\over 2^{(2-\delta)}(2\!-\!\delta)}
\, t^{\delta-3}
$$
and with $\delta=2(1\!+\!\nu)<0$ and $\mu=1\!-\!\nu$. 
However, these densities are in terms of new classes of special 
functions. The clarification of their relations to those for Bessel 
processes of non-negative indices is but one of the problems 
that require further study.
\vskip.9cm
\centerline{\Large\bf V\quad Laplace transforms in the general case:} 
\vskip.06cm
\centerline{\Large\bf combining stochastics and complex analysis}
\vskip.3cm\noindent
{\large\bf 15.\quad Remarks on general philosphy:}\quad
%
The discussion up to now has focused in particular on extension 
at the process level. The actual valuation problem, however,
is not at the process level but at the expectation level.
From this point of view, \S 10 has identified two function $f$ and
$g$ in the variable $\nu$ and has proved them to be equal for 
$\nu$ non-negative. One would like to have this equality also
for negative $\nu$, and thus extend the validity of the identity 
$f=g$ from the non-negative real line to the whole real line.     
Such situations quite commonly appear in problems 
in analysis and are addressed there using analytic continuation.
However, functions become amenable to complex analytic methods
only on open subsets of the complex plane. Thus the {\it identity 
theorem\/} of complex analysis asserts that two functions on a 
connected open subset 
of the complex plane are equal if they are analytic and agree on a 
convergent sequence there only. So it is in fact no longer possible 
to stick to real numbers only. As a subset of the complex plane they 
are closed with an empty interior. At this stage, however, nothing 
is known about the functions of \S 10 if $\nu$ is outside the 
non-negative real line. In particular, it is not known if they exist 
at all, and this needs to be established togother with their 
analyticity properties. Of these two functions, $D_\nu(a,z)$
is already given as an explicit analytic expression, while the other 
function is not. In fact, it is not explicit at all, as it is defined 
as the Laplace transform of a certain expectation. The question to 
be tackled here is then how to get explicit analyticity properties 
from such non-explicit  stochastic concepts. In the present situation, 
we attack this in two stages.
First, establish analyticity properties of the expectation. Then,
in a second step, study how these are preserved on taking Laplace
tranforms. As it turns out,  using standard results from complex 
analysis, we find that it is not Bessel precesses which
enter into the analysis but simply Brownian motion.  All of this 
may be regarded as an instance for why
enhancing stochastics by complex analytic 
methods seems quite an interesting and promising line of thought. 
\goodbreak\grosserAbsatz
%
%
{\large\bf 16.\quad First step -- analyticity of the function 
$D_\nu(a,z)$:}\quad
As a first step in the analytic continuation argument, this section 
studies the ana\-lytic properties of the generalized first Weber 
integral $D_\nu(a,z)$ of \S 8~Theorem. On choosing the square root 
associated to the principal branch of the logarithm recall 
$$
D_\nu(a,z)
=
{e\vbox to 6pt{}^{ -{\scriptstyle 1\over\scriptstyle 2a}}
\over a}
\int_0^\infty 
e\vbox to 8pt{}^{-{\scriptstyle x^2\over\scriptstyle 2a}}
x^{\nu+3}
I_{\sqrt{2z\!+\! \nu^2}}\Big( {x\over a}\Big)\, dx\, $$
is finite for any positive real number $a$ and for any complex numbers 
$z$ and $\nu$ such that the real part of $\nu\!+\!4\!+\!(2z\!+\!\nu^2)^{1/2}$ 
is positive as a consequence of \S 11~Proposition's integrability analysis. 
Using the confluent hypergeometric function $\Phi$ 
discussed in \cite[\S 9.9]{Lebedev}, the precise analyticity result to be 
proved is the following   
\grosserAbsatz
{\bf Proposition:}\quad{\it Let $a$ be any positive real and 
$\varepsilon$ be any non-negative real. For any complex number 
$\nu$ with $|\im(\nu)|\le \varepsilon$  we then have
$$
D_\nu(a,z)=
\Gamma\Big( {\nu\!+\!4\!+\!\mu\over 2}\Big)
\cdot{1\over \Gamma(\mu\!+\!1)}
  \Phi\Big( {\nu\!+\!4\!+\!\mu\over 2 }, \mu\!+\!1; {1\over 2a}\Big)
\cdot e\vbox to 6pt{}^{ -{\scriptstyle 1\over\scriptstyle 2a}}
\cdot           
(2a)\vbox to 9pt{}^{{\scriptstyle \nu+2-\mu\over\scriptstyle 2}}
$$
if $z$ is any complex with $\re(z)>2\varepsilon^2$ setting  
$\mu=\sqrt{2z\!+\!\nu^2\, }$.} 
\goodbreak\grosserAbsatz
{\bf Corollary:}\quad{\it Let $a$ be any positive real and 
$\varepsilon$ be any non-negative real. 
For any complex number $z$ with $\re(z)>2\varepsilon^2$, 
sending $\nu$ to $D_\nu(a,z)$ then defines
an analytic map on the set of all complex numbers $\nu$ with 
$|\im(\nu)|\le\varepsilon$.}
\goodbreak\mittlererAbsatz
Both results are based on \S 11~Proposition which gives finiteness 
of $D_\nu(a,z)$ under their conditions on $\nu$, $a$, and $z$. 
Combining this with the analyticity properties of $\Phi$ discussed 
in \cite[\S 9.9]{Lebedev} and those of the gamma function, the 
Corollary follows from the Proposition.  The proof of the Proposition
then reduces to exlicitly computing $D_\nu(a,z)$. For this 
we modify the quite typical discussion in 
\cite[\S13.3, pp.393f]{Watson} of Hankel's generalization of 
Weber's first integral. The idea is to expand the modified 
Bessel function in the integrand of 
$$
I=\int_0^\infty 
e\vbox to 9pt{}^{ -{\scriptstyle x^2\over \scriptstyle 2a}}
x^{\nu+3}I_\mu\Big( {x\over a}\Big) dx
$$
into its series of \cite[\S 5.7]{Lebedev} and integrate term by 
term. Using \cite[\S 9.9]{Lebedev} this is justified by the absolute 
convergence of the series for the confluent hypergeometric series 
which is to result, and we get
$$I
={1\over (2a)^\mu} 
\sum_{n=0}^\infty 
{1\over \Gamma(\mu\!+\!1\!+\!n)} {(2a)^{-2n} \over n!}
\int_0^\infty 
 e\vbox to 9pt{}^{ -{\scriptstyle x^2\over \scriptstyle 2a}}
x\vbox to 8pt{}^{ \nu+3+\mu+2n}\, dx\,. $$
\goodbreak\noindent
Changing variables $y=(2a)^{-1}x^2$, compute any $n$-th integral as 
$$
\int_0^\infty 
 e\vbox to 9pt{}^{ -{\scriptstyle x^2\over \scriptstyle 2a}}
x\vbox to 8pt{}^{ \nu+3+\mu+2n}\, dx\
=
{1\over 2}\cdot
\Gamma\Big( 
{\nu\!+\!\mu\!+\!4\over \scriptstyle 2}+n\Big)
\cdot(2a)\vbox to 8pt{}^{ {\scriptstyle \nu+\mu+4\over \scriptstyle 2}+n}. 
$$
Extracting the series of the pertinent confluent hypergeometric 
function we thus get 
$$I= {1\over 2}\cdot
{\Gamma\big( (\nu\!+\!4\!+\!\mu)/2\big)\over \Gamma(\mu\!+\!1)}
\cdot \Phi\Big( {\nu\!+\!4\!+\!\mu\over 2 }, \mu\!+\!1; {1\over 2a}\Big)
\cdot (2a)\vbox to 8pt{}^{ {\scriptstyle \nu-\mu+4\over \scriptstyle 2}} .$$
Multiplying this expression with 
$\exp( -( 2a)^{-1})/a$, the Proposition follows. 
%
\goodbreak\grosserAbsatz
%
{\large\bf 17.\quad Second step -- Analyticity of the 
functions $f_{GY,a}$:}\quad
Establishing analyticity results in $\nu$ about the Laplace 
transform of the expectation defining non-Asian option prices 
combines insights from stochastics with insights of an analytic
nature. We establish this result in two steps. As a first step, 
in this section any of the auxiliary functions $f_{GY,a}$
of \S 7 is considered as function in the variable $\nu$ 
$$  L(x,\nu)=E^{Q}\big[ \big( A^{(\nu)}_x-a\big)^{\!+}\big] 
$$ 
for any fixed positive real numbers $a$ and $x$. Using \S 10~Lemma, 
we know it is defined for non-negative real numbers $\nu$. 
However, this has been achieved in a very indirect way only: for 
these values of $\nu$ the Laplace transforms of the corresponding 
functions $f_{GY,a}$ have been shown to be finite. Now more is true 
indeed
\grosserAbsatz
{\bf Lemma:}\quad{\it For any $x>0$, the function $\nu\mapsto L(x,\nu)$ 
extends to a function on the complex plane which is analytic at each 
point, and  for which we have the majorization}
$$ |L(x,\nu)|\le e\vbox to 8pt{}^{
{\scriptstyle x\over\scriptstyle  2}\im^2(\nu)}
             \cdot E^Q[ A_x^{(\re(\nu))}].$$
\kleinerAbsatz
For the proof of the Lemma now set 
$f(w)=(w\!-\!a)^+$. Applying Girsanov's theorem such that 
$W_x=\nu x\!+\!B_x$ becomes a standard Brownian motion, and 
dropping reference to this new measure, we get
$$ L(x,\nu)= E\Big[f\big(A_x^{(0)}\big)e\vbox to 8pt{}^{\nu W_x}\Big]
\cdot e\vbox to 8pt{}^{-{\scriptstyle x\over\scriptstyle  2}\nu^2} 
.$$
For establishing the analyticity statement of the Lemma, it is thus
sufficient to show that the expectation factor is analytic in any 
complex number $\nu$. This is true by definition if we have the 
convergent series 
$$
E\Big[f\big(A_x^{(0)}\big)e\vbox to 8pt{}^{\nu W_x}\Big]
=
\sum_{m=0}^\infty 
{\nu^m\over m!}E\Big[f\big(A_x^{(0)}\big) W_x^m\Big]
$$
for all $\nu$. For this it is sufficient to show that the series 
is absolutely convergent for all $\nu$. Using the Cauchy-Schwarz 
inequality this is implied by the convergence of 
$$
\sum_{m=0}^\infty 
{|\nu|^m\over m!}\sqrt{E\Big[f^2\big(A_x^{(0)}\big)\Big]\, }
\sqrt{E\Big[W_x^{2m}\Big]\, }
$$
for all $\nu$. Herein $E[f^2(A_x^{\raise-2pt\hbox{$\scriptstyle (0)$}})]$
is majorized by the second moment of Yor's zero drift  process 
$A^{\raise-2pt\hbox{$\scriptstyle (0)$}}$ at $x$, and so is finite 
from \cite[\S4]{Yor 92}. Since  the factors $E[W_x^{2m}]$ are majorized by  
$\pi^{-1/2}\cdot (2x)^m\cdot m!$ for all $m\ge 0$, convergence follows 
using the ratio test. 
Actually, we have so established yet another upper bound
to the price of the Asian option. 
\mittlererAbsatz
To establish the majorization of the Lemma, taking
absolute values inside the expectation in the above Girsanov
representation of $L(x,\nu)$ gives:
$$|L(x,\nu)|\le 
E\Big[f\big(A_x^{(0)}\big)\cdot \big|e\vbox to 8pt{}^{\nu W_x}\big|\, \Big]
\cdot \big|e\vbox to 8pt{}^{-{\scriptstyle x\over\scriptstyle  2}\nu^2}\big| 
.$$
The absolute value of the exponential factors are the exponentials of 
the real parts of the respective arguments.
Majorizing the function in Yor's zero drift process 
by this process itself, we get
$$|L(x,\nu)|\le 
e\vbox to 8pt{}^{{\scriptstyle x\over\scriptstyle  2}\im^2(\nu)}\cdot 
E\Big[A_x^{(0)}e\vbox to 8pt{}^{\re(\nu) W_x}\Big]
e\vbox to 8pt{}^{-{\scriptstyle x\over\scriptstyle  2}\re^2(\nu)}
.$$ 
Reversing the Girsanov transformation then completes the proof of the Lemma.
\goodbreak\grosserAbsatz
%
{\large\bf 18.\quad Third step -- analyticity of the transforms $F_{GY,a}$:}\quad
While the preceding section has studied the expectations 
$$ 
L(x,\nu)=E^{Q}\big[ \big( A^{(\nu)}_x-a\big)^{\!+}\big] 
$$ 
for any fixed positive real numbers $a$ and $x$ as function in the 
complex variable $\nu$ only, this section moreover treats $x$ as a 
variable. If $\nu$ is any non-negative real number, we have from 
\S 10~Lemma that the integrals of the Laplace transforms
$$F(\nu)(z)=\int_0^\infty e^{-zx}L(x,\nu)\, dx
$$ 
are finite if $\re(z)>2(\nu\!+\!1)$. In this section we give 
an independent proof of the following more general statement
\grosserAbsatz
{\bf Proposition:}\quad{\it For any complex number $z$
with a positive real part, the map sending $\nu$ to $F(\nu)(z)$ 
is analytic in all complex numbers $\nu$ with  
$\re(z)>{1\over 2}\im^2(\nu)\!+\!2(\re(\nu)\!+\!1)$.} 
\grosserAbsatz
The proof of the Proposition is based on the following
\grosserAbsatz
{\bf Lemma:}\quad{\it For any complex number $\nu$, the Laplace 
transform $F(\nu)(z)$ is finite for any complex number $z$ with 
$\re(z)>\max\{ 0,{1\over 2}\im^2(\nu)\!+\!2(\re(\nu)\!+\!1)\}$.} 
\mittlererAbsatz
{\bf Proof of the Lemma:}\quad As first step in proving the 
Lemma, we establish for any complex number $\nu$ the majorization
$$
|F(\nu)(z)|\le \int_0^\infty e\vbox to 8pt{}^{
             -\big(\re(z)-{\scriptstyle 1\over\scriptstyle 2}\im^2(\nu)\big)x} 
E^Q[A_x^{(\re(\nu))}]\, dx$$
for any complex number $z$ in the sense of measurable functions. 
Indeed, majorize the absolute value of  $F(\nu)(z)$ by taking 
the absolute value inside the defining integral. The absolute value of 
the exponential factor then is equal to $\exp(-\re(z)x)$. 
Majorizing the absolute value of  $L(x,\nu)$ using \S 17~Lemma,
the estimate follows. 
\kleinerAbsatz
Setting $\nu_0=\re(\nu)$, now let $\re(z)$ be positive and such that
$\xi_0=\re(z)\!-\!\im(\nu)^2/2$ is bigger than $2(\nu_0\!+\!1)$. We 
compute the Laplace transform of the right hand side of the above 
inequality using \S 4~Lemma.
If $\nu_0$ is different from minus one, 
$$
\int_0^\infty e^{-\xi_0x} E^Q[A_x^{(\nu_0)}]\, dx
=
{1\over \xi_0(\xi_0-2(\nu_0\!+\!1))}\, , 
$$
using that  $\xi_0$ is bigger than $2(\nu_0\!+\!1)$ to compute the 
improper integrals. It converges to $\xi_0^{-2}$ with $\nu_0$ going 
to minus one. Thus it is seen to coincide with the Laplace transform 
for the case $\nu_0=-1$. 
The Laplace transforms are finite if $\xi_0$ is positive 
and bigger than $2(\nu_0\!+\!1)$, and then $F(\nu)(z)$ is finite
{\it a forteriori\/}. The proof of the Lemma is complete.
\grosserAbsatz
{\bf Proof of the Proposition:}\quad For proving the Proposition,
fix any complex $z$ with a positive real part and then choose any 
complex $\nu_0$ satisfying the resulting inequality of the 
Proposition. In fact, the validity of this inequality then 
extends to all $\nu$ in a compact neighbourhood $V$ of $\nu_0$. 
As a consequence of \S 17~Lemma, $f_z(x,\nu)=\exp(-zx)L(x,\nu)$ 
is for any $x>0$ analytic in $\nu$ on $V$ in particular. 
Since $V$ is compact, the argument of the Lemma moreover shows 
that its absolute value is majorized by an integrable function $g$ 
on the real line. The function $\nu\mapsto F(\nu)(z)$ so is 
continuous on $V$ being obtained by integration of $f_z$ over 
the positive real line. To show it is analytic on the interior 
of $V$ we want to apply {\it Morera's theorem\/}, see 
\cite[10.17, p.208]{Rudin} and have to show
$$
\int_{\partial \Delta} F(\nu)(z)\, d\nu =0
$$
for any triangle $\Delta$ contained in the interior of $V$. 
Indeed, since we have shown $f_z$ to be integrable, applying Fubini's 
theorem gives 
$$
\int_{\partial \Delta} F(\nu)(z)\, d\nu
=
\int_0^\infty\int_{\partial \Delta}
f_z(x,\nu)\,d\nu\, dx\, .
$$
Recalling that $\nu\mapsto f_z(x,\nu)$ is analytic from \S 17~Lemma, 
the inner integral herein is zero by {\it Cauchy's theorem\/}
for $\Delta$, see \cite[10.13, p.205]{Rudin}. Thus the whole 
double integral is zero, as was to  be shown. This completes 
the proof of the Proposition         
\goodbreak\grosserAbsatz
%
{\large\bf 19.\quad Final step -- second proof of the Laplace transform using 
analytic extension:}\quad 
It remains to pull things together and show how the results of this 
Part V give a second proof of the two results of \S 8 by using 
analytic extension.
\mittlererAbsatz
First notice that \S 8~Corollary is implied by \S 8~Theorem
using the computation of $D_\nu(a,z)$ in \S 16~Lemma. Thus we
are reduced to prove \S 8~Theorem.
\mittlererAbsatz
The proof of \S 8~Theorem is by analytic continuation using the 
identity theorem, see \cite[Corollary to 10.18, p.209]{Rudin}. 
As a consequence of \S 10~Lemma 
it remains to establish the identity in \S 8~Theorem for negative 
indices $\nu$ only. The existence of the Laplace transform on all 
complex numbers $z$ with positive real part required in the Theorem 
then has essentially been proved in \S 18~Proposition. 
\mittlererAbsatz
For establishing the crucial Laplace transform identity of \S 8~Theorem,
thus first let $z$ be any complex number with real part $\re(z)>4$ and choose
$0<\varepsilon<1$. Using \S 16~Corollary, the generalized first Weber 
integral
$$
D_\nu(a,z)
=
{e\vbox to 6pt{}^{ -{\scriptstyle 1\over\scriptstyle 2a}}
\over a}
\int_0^\infty 
e\vbox to 8pt{}^{-{\scriptstyle x^2\over\scriptstyle 2a}}
x^{\nu+3}
I_{\sqrt{2z\!+\! \nu^2}}\Big( {x\over a}\Big)dx\, 
$$ 
is analytic in $\nu$ on the $\varepsilon$-thickened real line 
$A_\varepsilon$ which consists of all complex numbers $\nu$ with
$|\im(\nu)|<\varepsilon$. Picture $A_\varepsilon$ as the band of 
height $2\varepsilon$ symmetric with respect to the real axis.
If we choose $\varepsilon$ so small that $\re(z)>2(2\!+\!2\varepsilon)$, we
claim to have analyticity of the Laplace transform 
$$
F(\nu)(z)=\int_0^\infty e^{-zx}L(x,\nu)\, dx
$$     
as a function in $\nu$ on the $\varepsilon$-thickened half-line 
$B_\varepsilon$, that subset of $A_\varepsilon$
which consists of  
\[
\xy
/r0.4cm/:
 (0,0)*{\bullet}="0",
 (-.4,-.4)*+!DR{0}, 
 (1,1)*{\bullet}="u",
 (1,-1)*{\bullet}="d",
 (1.1,1.1)*!DL{(\varepsilon,\varepsilon)},
 (1.1,-1.1)*!UL{(\varepsilon,-\varepsilon)},
 \ar @*{[|(2)]}  (-6,0);(3,0),
 \ar @*{[|(2)]}  (0,-3);(0,3)
 \ar@{-}@*{[|(3)]}  (-6,1);"u",
 \ar@{-}@*{[|(3)]}  (-6,-1);"d",
 \ar@{-}@*{[|(3)]}  "u";"d",
\ar@{-} (-6.000000000000000000000000000  ,  -1  )  ;  (-6.000000000000000000000000000   ,  1  ),
\ar@{-} (-5.900000000000000000000000000  ,  -1  )  ;  (-5.900000000000000000000000000   ,  1  ),
\ar@{-} (-5.799999999999999999999999999  ,  -1  )  ;  (-5.799999999999999999999999999   ,  1  ),
\ar@{-} (-5.699999999999999999999999999  ,  -1  )  ;  (-5.699999999999999999999999999   ,  1  ),
\ar@{-} (-5.599999999999999999999999999  ,  -1  )  ;  (-5.599999999999999999999999999   ,  1  ),
\ar@{-} (-5.500000000000000000000000000  ,  -1  )  ;  (-5.500000000000000000000000000   ,  1  ),
\ar@{-} (-5.400000000000000000000000000  ,  -1  )  ;  (-5.400000000000000000000000000   ,  1  ),
\ar@{-} (-5.299999999999999999999999999  ,  -1  )  ;  (-5.299999999999999999999999999   ,  1  ),
\ar@{-} (-5.199999999999999999999999999  ,  -1  )  ;  (-5.199999999999999999999999999   ,  1  ),
\ar@{-} (-5.099999999999999999999999999  ,  -1  )  ;  (-5.099999999999999999999999999   ,  1  ),
\ar@{-} (-5.000000000000000000000000000  ,  -1  )  ;  (-5.000000000000000000000000000   ,  1  ),
\ar@{-} (-4.900000000000000000000000000  ,  -1  )  ;  (-4.900000000000000000000000000   ,  1  ),
\ar@{-} (-4.800000000000000000000000000  ,  -1  )  ;  (-4.800000000000000000000000000   ,  1  ),
\ar@{-} (-4.699999999999999999999999999  ,  -1  )  ;  (-4.699999999999999999999999999   ,  1  ),
\ar@{-} (-4.599999999999999999999999999  ,  -1  )  ;  (-4.599999999999999999999999999   ,  1  ),
\ar@{-} (-4.500000000000000000000000000  ,  -1  )  ;  (-4.500000000000000000000000000   ,  1  ),
\ar@{-} (-4.400000000000000000000000000  ,  -1  )  ;  (-4.400000000000000000000000000   ,  1  ),
\ar@{-} (-4.300000000000000000000000000  ,  -1  )  ;  (-4.300000000000000000000000000   ,  1  ),
\ar@{-} (-4.199999999999999999999999999  ,  -1  )  ;  (-4.199999999999999999999999999   ,  1  ),
\ar@{-} (-4.099999999999999999999999999  ,  -1  )  ;  (-4.099999999999999999999999999   ,  1  ),
\ar@{-} (-4.000000000000000000000000000  ,  -1  )  ;  (-4.000000000000000000000000000   ,  1  ),
\ar@{-} (-3.899999999999999999999999999  ,  -1  )  ;  (-3.899999999999999999999999999   ,  1  ),
\ar@{-} (-3.799999999999999999999999999  ,  -1  )  ;  (-3.799999999999999999999999999   ,  1  ),
\ar@{-} (-3.699999999999999999999999999  ,  -1  )  ;  (-3.699999999999999999999999999   ,  1  ),
\ar@{-} (-3.599999999999999999999999999  ,  -1  )  ;  (-3.599999999999999999999999999   ,  1  ),
\ar@{-} (-3.500000000000000000000000000  ,  -1  )  ;  (-3.500000000000000000000000000   ,  1  ),
\ar@{-} (-3.399999999999999999999999999  ,  -1  )  ;  (-3.399999999999999999999999999   ,  1  ),
\ar@{-} (-3.299999999999999999999999999  ,  -1  )  ;  (-3.299999999999999999999999999   ,  1  ),
\ar@{-} (-3.200000000000000000000000000  ,  -1  )  ;  (-3.200000000000000000000000000   ,  1  ),
\ar@{-} (-3.099999999999999999999999999  ,  -1  )  ;  (-3.099999999999999999999999999   ,  1  ),
\ar@{-} (-3.000000000000000000000000000  ,  -1  )  ;  (-3.000000000000000000000000000   ,  1  ),
\ar@{-} (-2.899999999999999999999999999  ,  -1  )  ;  (-2.899999999999999999999999999   ,  1  ),
\ar@{-} (-2.799999999999999999999999999  ,  -1  )  ;  (-2.799999999999999999999999999   ,  1  ),
\ar@{-} (-2.700000000000000000000000000  ,  -1  )  ;  (-2.700000000000000000000000000   ,  1  ),
\ar@{-} (-2.599999999999999999999999999  ,  -1  )  ;  (-2.599999999999999999999999999   ,  1  ),
\ar@{-} (-2.500000000000000000000000000  ,  -1  )  ;  (-2.500000000000000000000000000   ,  1  ),
\ar@{-} (-2.400000000000000000000000000  ,  -1  )  ;  (-2.400000000000000000000000000   ,  1  ),
\ar@{-} (-2.299999999999999999999999999  ,  -1  )  ;  (-2.299999999999999999999999999   ,  1  ),
\ar@{-} (-2.200000000000000000000000000  ,  -1  )  ;  (-2.200000000000000000000000000   ,  1  ),
\ar@{-} (-2.099999999999999999999999999  ,  -1  )  ;  (-2.099999999999999999999999999   ,  1  ),
\ar@{-} (-2.000000000000000000000000000  ,  -1  )  ;  (-2.000000000000000000000000000   ,  1  ),
\ar@{-} (-1.899999999999999999999999999  ,  -1  )  ;  (-1.899999999999999999999999999   ,  1  ),
\ar@{-} (-1.799999999999999999999999999  ,  -1  )  ;  (-1.799999999999999999999999999   ,  1  ),
\ar@{-} (-1.699999999999999999999999999  ,  -1  )  ;  (-1.699999999999999999999999999   ,  1  ),
\ar@{-} (-1.600000000000000000000000000  ,  -1  )  ;  (-1.600000000000000000000000000   ,  1  ),
\ar@{-} (-1.500000000000000000000000000  ,  -1  )  ;  (-1.500000000000000000000000000   ,  1  ),
\ar@{-} (-1.399999999999999999999999999  ,  -1  )  ;  (-1.399999999999999999999999999   ,  1  ),
\ar@{-} (-1.299999999999999999999999999  ,  -1  )  ;  (-1.299999999999999999999999999   ,  1  ),
\ar@{-} (-1.200000000000000000000000000  ,  -1  )  ;  (-1.200000000000000000000000000   ,  1  ),
\ar@{-} (-1.100000000000000000000000000  ,  -1  )  ;  (-1.100000000000000000000000000   ,  1  ),
\ar@{-} (-1.000000000000000000000000000  ,  -1  )  ;  (-1.000000000000000000000000000   ,  1  ),
\ar@{-} (-0.8999999999999999999999999999  ,  -1  )  ;  (-0.8999999999999999999999999999   ,  1  ),
\ar@{-} (-0.8000000000000000000000000000  ,  -1  )  ;  (-0.8000000000000000000000000000   ,  1  ),
\ar@{-} (-0.6999999999999999999999999999  ,  -1  )  ;  (-0.6999999999999999999999999999   ,  1  ),
\ar@{-} (-0.6000000000000000000000000000  ,  -1  )  ;  (-0.6000000000000000000000000000   ,  1  ),
\ar@{-} (-0.5000000000000000000000000000  ,  -1  )  ;  (-0.5000000000000000000000000000   ,  1  ),
\ar@{-} (-0.4000000000000000000000000000  ,  -1  )  ;  (-0.4000000000000000000000000000   ,  1  ),
\ar@{-} (-0.3000000000000000000000000000  ,  -1  )  ;  (-0.3000000000000000000000000000   ,  1  ),
\ar@{-} (-0.2000000000000000000000000000  ,  -1  )  ;  (-0.2000000000000000000000000000   ,  1  ),
\ar@{-} (-0.1000000000000000000000000000  ,  -1  )  ;  (-0.1000000000000000000000000000   ,  1  ),
\ar@{-} (0  ,  -1  )  ;  (0   ,  1  ),
\ar@{-} (0.09999999999999999999999999999  ,  -1  )  ;  (0.09999999999999999999999999999   ,  1  ),
\ar@{-} (0.1999999999999999999999999999  ,  -1  )  ;  (0.1999999999999999999999999999   ,  1  ),
\ar@{-} (0.2999999999999999999999999999  ,  -1  )  ;  (0.2999999999999999999999999999   ,  1  ),
\ar@{-} (0.3999999999999999999999999999  ,  -1  )  ;  (0.3999999999999999999999999999   ,  1  ),
\ar@{-} (0.4999999999999999999999999999  ,  -1  )  ;  (0.4999999999999999999999999999   ,  1  ),
\ar@{-} (0.5999999999999999999999999999  ,  -1  )  ;  (0.5999999999999999999999999999   ,  1  ),
\ar@{-} (0.6999999999999999999999999999  ,  -1  )  ;  (0.6999999999999999999999999999   ,  1  ),
\ar@{-} (0.7999999999999999999999999999  ,  -1  )  ;  (0.7999999999999999999999999999   ,  1  ),
\ar@{-} (0.8999999999999999999999999999  ,  -1  )  ;  (0.8999999999999999999999999999   ,  1  ),
\ar@{-} (0.9999999999999999999999999999  ,  -1  )  ;  (0.9999999999999999999999999999   ,  1  ),
 \endxy
\]
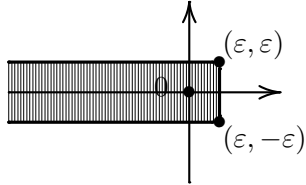
\centerline{{\ninerm}{\ninebf Figure 1.}\enspace The 
${\varepsilon}$-thickened half-line $\scriptstyle B_\varepsilon$.}
\vskip10pt\noindent
all complex numbers $\nu$ with $|\im(\nu)|<\varepsilon$
and $\re(\nu)<\varepsilon$.
Indeed, if $\re(\nu)<\varepsilon$, we have 
$2(2\!+\!2\varepsilon)>2\varepsilon\!+\! 2(\re(\nu)\!+\!2)$. If 
$|\im(\nu)|<\varepsilon$, we have $2\varepsilon>\im^2(\nu)/2$ since 
$\varepsilon<1$. Any $\nu$ in $B_\varepsilon$ so satisfies the 
inequality $\re(z)> \im^2(\nu)/2\!+\!2(\re(\nu)\!+\!2)$ of 
\S 18~Proposition, and the claim follows. 
\mittlererAbsatz
So the functions $F(\nu)(z)$ and 
$$ D_\nu^*(a,z)
=
{D_\nu(a,z)\over z(z-2(\nu\!+\!1))}
$$
are analytic as functions in $\nu$ on the $\varepsilon$-thickened 
half-line $B_\varepsilon$.
It is now a consequence of \S 10~Lemma  that we have  
$F(\nu)(z)=D_\nu^*(a,z)$ for all $\nu\ge0$ such that
$2(\nu\!+\!1)<\re(z)$. With $\re(z)>4$ this so holds 
{\it a forteriori\/} for all $\nu\ge 0$ such that 
$2(\nu\!+\!1)<4$, i.e., for all non-negative real numbers 
$\nu$ smaller than $1$. With $\varepsilon$ smaller than
$1$  we so have  $F(\nu)(z)=D_\nu^*(a,z)$ 
for any $\nu$ in the subinterval $(0,\varepsilon)$ of 
$B_\varepsilon$. With $B_\varepsilon$ open and connected 
the identity theorem literally applies to give that 
$F(\nu)(z)=D^*_\nu(a,z)$ for all $\nu$ in $B_\varepsilon$. 
This identity then holds {\it a forteriori\/} for all real 
numbers $\nu$ in $B_\varepsilon$, i.e., for all $\nu<\varepsilon$. 
This gives \S 8~Theorem for any $z$ with $\re(z)>4$.
\mittlererAbsatz
To lift this last restriction on $\re(z)$, notice that 
\S 18~Lemma implies $F(\nu)$ for any fixed real $\nu$ 
to be analytic on the half-plane $\{\re(z)>2(\nu\!+\!1)\}$;
this is seen by using a Morera type argument as for proving
\S 18~Proposition.  
On the other hand, applying \S 16~Corollary, $D_\nu^*(a,z)$ 
as a function of $z$ is analytic on the intersection of this 
last half-plane with the right-half plane. The validity of 
the identity 
$F(\nu)(z)=D_\nu^*(a,z)$ 
thus can be analytically continued from complex numbers $z$ 
with $\re(z)>4$ to complex numbers $z$ with $\re(z)$ positive 
and bigger than $2(\nu\!+\!1)$, and the second proof of 
\S 8~Theorem is complete. 
\vskip.9cm
\centerline{\Large\bf VI\quad Laplace transforms in the general case:} 
\centerline{\Large\bf a uniform proof}
\vskip.3cm\noindent
%
{\large\bf 20.\quad Remarks on general philosophy:}\quad
The question may arise if there is not a uniform way for computing 
the Laplace transfoms in \S 8~Theorem not requiring a two step 
procedure. The Lamperti factorization  
$$
e\vbox to 8pt{}^{B_w+\nu w}=R^{(\nu)}\big( A^{(\nu)}(w)\big)
$$
of \S 9~Lemma may in fact offer a clue. The idea of the 
two arguments discussed up to now is to focus on the Bessel 
process side of this identity. It is the problems with 
Bessel processes that transcribe into their two-step 
approach. A possible remedy so would be to try to bring 
in Bessel processes not at the earliest stage but as late 
as possible. The Lamperti factorization suggests to focus 
on the geometric Brownian motion of its left hand side.      
While Yor has already made extensive use of this device for 
instance in [{\bf Y92}, \S6], we explain in the sequel what seems 
to be a most natural adaption of the Girsanov technique to the 
Lamperti factorization approach for explicitly computing 
the Laplace transforms $F_{GY,a}$ with $a>0$. 
%
%
\goodbreak\grosserAbsatz
{\large\bf 21.\quad Third proof of the theorem -- a uniform argument:}\quad
This section sketches a uniform way for computing the Laplace 
transforms $F_{GY,a}$ with $a$ any positive real, and thus 
provides a third proof of \S 8~Theorem. 
Recall these transforms are defined by
$$F_{GY,a}(z)
=\int_0^\infty e\vbox to 8pt{}^{ -zx} f_{GY,a}(x)\, dx
$$
for any complex $z$ with $\re(z)$ sufficiently big, 
and where
$$
f_{GY, a}(x)=E^Q\big[ f\big(A^{\raise-1.4pt\hbox{$\scriptstyle(\nu)$}}(x)\big)\big].
$$
setting $f(x)=(x\!-\!a)^+$. The key idea is to apply 
Girsanov's theorem such that $W_x=\nu x\!+\!B_x$ becomes 
a Brownian motion, and suppressing reference to this new measure, 
transcribe $f_{GY,a}$ in terms of Yor's zero drift process 
$A^{\raise-1.4pt\hbox{$\scriptstyle(0)$}}$ as follows 
$$
f_{GY,a}(x)
=
E^Q\Big[f\big(A^{\raise-1.4pt\hbox{$\scriptstyle(0)$}}_x\big)
e\vbox to 9pt{}^{ -{\nu^2\over 2}x+\nu W_x}\Big].
$$ 
With 
$\tau_{0,a}=\inf\{u\, |\, A^{\raise-1.4pt\hbox{$\scriptstyle(0)$}}(u)>0\}$ 
the inverse time change to $A^{\raise-1.4pt\hbox{$\scriptstyle(0)$}}$ at 
time zero, the computations of \S 10
thus prove $f_{GY,a}(x)=0$ on the set of all events where $x\le\tau_{0,a}$. 
And on the set of all events where $x\ge \tau_{0,a}$ we now have the 
representation 
$$f\big(A^{\raise-1.4pt\hbox{$\scriptstyle(0)$}}\big)
= \big(R^{\raise-1.4pt\hbox{$\scriptstyle(0)$}}_a\big)^2
\cdot A^{\raise-1.4pt\hbox{$\scriptstyle(0)$}}_{x\!-\!\tau_{0,a}}, 
$$
which applies to express $f_{GY,a}(x)$ as the following 
iterated expecation
$$
f_{GY,a}(x)=E\bigg[
 \big(R^{\raise-1.4pt\hbox{$\scriptstyle(0)$}}_a\big)^2 
E\bigg[  
A^{\raise-1.4pt\hbox{$\scriptstyle(0)$}}_{x\!-\!\tau_{0,a}}  
e\vbox to 9pt{}^{ -{\nu^2\over 2}x+\nu W_x}
\Big|\,
 \Fgeschwungen _{\tau_{0,a}}
\bigg]
\bigg]\, . 
$$
Applying Laplace transforms to both sides of this identity
then gives
$$
F_{GY,a}(z)
=
E\bigg[ 
\big(R^{\raise-1.4pt\hbox{$\scriptstyle(0)$}}_a\big)^2
e\vbox to 9pt{}^{ -\big(z+{\nu^2\over 2}\big)\tau_{0,a}}
E\bigg[ \int_0^\infty e\vbox to 8pt{}^{-zx} 
 A^{\raise-1.4pt\hbox{$\scriptstyle(0)$}}_{x}  
e\vbox to 9pt{}^{-{\nu^2\over 2}x+\nu W_{x+\tau_{0,a}}}\, dx
\Big|\,
 \Fgeschwungen _{\tau_{0,a}}\bigg]
\bigg]\, . 
$$
In the conditional expectation restart the Brownian motion 
at time $\tau_{0,a}$. To the resulting additional $\nu$-th 
power of $\exp(W_{\tau_0,a})$ apply the corresponding Lamperti 
factorization and interpret it as being equal to the 
$\nu$-th power of the zero index Bessel process at time $a$.
Collecting powers of this last process, 
$$
F_{GY,a}(z)  
=
E\bigg[ \big(R^{\raise-1.4pt\hbox{$\scriptstyle(0)$}}_a\big)^{\nu+2}
e\vbox to 9pt{}^{ -\big(z+{\nu^2\over 2}\big)\tau_{0,a}}
E\bigg[ \int_0^\infty e\vbox to 8pt{}^{-zx} 
 A^{\raise-1.4pt\hbox{$\scriptstyle(0)$}}_{x}  
e\vbox to 9pt{}^{-{\nu^2\over 2}x+\nu W_{x}}\, dx
\Big|\,
 \Fgeschwungen _{\tau_{0,a}}\bigg]
\bigg]\, . 
$$
The Strong Markov condition which we have used here also makes the whole 
integrand in the conditional expectation independent of time-$\tau_{0,a}$
information.  Reversing the Girsanov transformation in this 
inner expectation,
$$
F_{GY,a}(z)  
=
 E\bigg[ \big(R^{\raise-1.4pt\hbox{$\scriptstyle(0)$}}_a\big)^{\nu+2}
e\vbox to 9pt{}^{ -\big(z+{\nu^2\over 2}\big)\tau_{0,a}}
E\bigg[ \int_0^\infty e\vbox to 8pt{}^{-\big(z+{\nu^2\over 2}\big)x} 
 A^{\raise-1.4pt\hbox{$\scriptstyle(\nu)$}}_{x} \, dx
\bigg]
\bigg]\, .
$$
This puts us into the situation of \S 10. 
Partially reversing the Tonelli argument, the inner expectation 
is equal to the first moment
of the drift-$\nu$ process $A^{\raise-1.4pt\hbox{$\scriptstyle(\nu)$}}$ 
at time $x$ recalled in \S 4~Lemma. Computing the Laplace transform,
we thus arrive at
$$
F_{GY,a}(z)=
{1\over 2(\nu\!+\!1)}\bigg(
{1\over 1-2(\nu\!+\!1)}-{1\over z}\bigg)
E\bigg[ \big(R^{\raise-1.4pt\hbox{$\scriptstyle(0)$}}_a\big)^{\nu+2}
e\vbox to 9pt{}^{ -(z+{\nu^2\over 2})\tau_{0,a}}\bigg]\,
$$
if $\re(z)$ is bigger than at least $\max\{0,2(\nu\!+\!1)\}$. 
In this way we are reduced to compute the expectation factor.
This, however, proceeds as in \S 10 by conditioning on 
the corresponding Bessel process. The whole 
point is that this time conditioning is not on a Bessel 
process of an arbitrary index but on a Bessel processes 
of index zero. Observing
how $\nu$ enters now via the exponent of the index zero Bessel
process and via a shifting factor for the time change  
$\tau_{0,a}$, we so arrive at
$$
E\bigg[ \big(R^{\raise-1.4pt\hbox{$\scriptstyle(0)$}}_a\big)^{\nu+2}
e\vbox to 9pt{}^{ -\big(z+{\nu^2\over 2}\big)\tau_{0,a}}\bigg]
=
{e\vbox to 8pt{}^{ -{\scriptstyle 1\over\scriptstyle 2a}}\over a}
\int_0^\infty 
e\vbox to 8pt{}^{-{\scriptstyle x^2\over\scriptstyle 2a}}
x^{\nu+3}
I_{\sqrt{2z\!+\! \nu^2}}\Big( {x\over a}\Big)dx\, , 
$$
for any $z$ with sufficiently big positive real part. Using
\S 11~Proposition, this integral is finite if $\re(z)$ is 
positive and bigger than $2(\nu\!+\!1)$. 
So this third proof of \S 8~Theorem is complete.  
\vskip.9cm
\centerline{\Large\bf VII\quad Epilogue} 
\vskip.3cm\noindent
{\large\bf 22.\quad Consequences for Asian options -- Hermite functions:}\quad
Having persevered to this point, the reader may wonder about
the nature and the quality of the implications of the mathematics
developed up to now. Indeed, going back to the starting point 
of the journey, it seems that the Laplace transform approach makes
possible significant improvements in understanding \S 3's normalized 
prices $C^{(\nu)}(h,q)$ of Asian options themselves. These improvements
are on a structural level and as such make possible advances on
computing $C^{(\nu)}$ as a consequence. This is essentially by 
being able to establish new links of Asian option valuation 
with a well-studied class of special functions: the Hermite functions
to be reviewed in this section.       
\mittlererAbsatz
Following [{\bf L}, \S \S 10.2ff], to which we refer for details, 
{\it Hermite functions\/} $H_\mu$ are analytic on the complex plane 
as functions of both their variable $z$ and their degree $\mu$. 
If the real part $\re(\mu)$ of  $\mu$ is bigger than $-1$, they have 
the integral representation 
$$
H_\mu(z)
=
{2^{\mu+1}\over \sqrt{\pi\, }}e^{z^2}\int_0^\infty
e^{-x^2} x^{\mu}\cos\Big( 2zx\!-\!{1\over 2}\mu\pi\Big)\,dx\,  
.$$
So they specialize to the $\mu$-th Hermite polynomials if $\mu$ is 
any non-negative integer, whence
$H_0=1$, $H_1(z)=2z$, $H_2(z)=4z^2\!-\!2$, $H_3(z)=8z^3\!-\! 12 z$,
for example. If the real part of  $\mu$ is negative, however,
Hermite functions change their character. Then they have  the 
integral representation:
$$H_\mu(z)={1\over \Gamma(-\mu)}\int_0^\infty e^{-u^2-2zu}u^{-(\mu+1)}\, du \,
 ,$$
and so specialize via $(2/\sqrt{\pi})H_{-1}(z)=\exp(z^2)\erfc (z)$ 
to the complementary error function $\erfc$ recalled to be given by
$$
\erfc(z)={2\over \sqrt{\pi}}\int_z^\infty e\vbox to 9pt{}^{-\xi^2}d\xi
\, .$$ 
For any complex $\mu$, Hermite functions can be expressed in terms of 
the Kummer confluent hypergeometric function $\Phi$ by 
$$
H_\mu(z)={2^\mu\cdot\Gamma(1/2)\over \Gamma\big( (1\!-\!\mu)/2\big)}
\cdot\Phi\Big( -{\mu\over 2}, {1\over 2}; z^2\Big)
+
z\cdot {2^\mu\cdot \Gamma(-1/2)\over \Gamma\big( -\mu/2\big)}
\cdot\Phi\Big( {1\!-\!\mu\over 2}, {3\over 2}; z^2\Big)
$$
for any complex $z$. From this representation one can moreover see
how Hermite functions are connected with the parabolic cylinder 
functions $D_\mu$  and with the Kummer confluent hypergeometric 
function of the second kind $\Psi$. 
\goodbreak\grosserAbsatz
{\large\bf 23.\quad Consequences for Asian options -- new integral 
representations:}\quad
Hermite functions as recalled in the previous section 
appear naturally in the closed form 
solution for \S 3's normalized price $ C^{(\nu)}(h,q)$ of 
the Asian option we have developed in [{\bf SA}]. If the 
normalized strike price $q$ is positive, it expresses this 
value as the sum of integral representations which have a 
product structure. They are obtained by integrating 
products of Hermite functions $H_\mu$ with weighted 
error functions as follows: 
\grosserAbsatz
{\bf Theorem:}\quad {\it If  $q$ is positive, the normalized price 
$ C^{(\nu)}(h,q)$ of the Asian option is given by the following
difference 
$$
C^{(\nu)}(h,q)=ce\vbox to 7pt{}^{ 2h(\nu+1)}\,S_{\nu+2}- c\,S_\nu\,  
$$
where the $S_\xi$ are three-term sums
$$
S_\xi= \ctrig {\xi} (\rho_\xi) +\chyp {\xi} (\rho_\xi)
+\chyp {-\xi} (\rho_\xi)
$$
whose single summands are integrals that depend on parameters
$\rho_\xi\ge 0$, but which as a whole are independent of these.} 
\grosserAbsatz
In terms of \S 3's concepts, $c$ is given by
$$
c=c(\nu,q)=
{\Gamma(\nu\!+\!4)\cdot 
(2q)\vbox to 7pt{}^{{1\over 2}{\scriptstyle(\nu+2)}}  
\over 2\pi \cdot (\nu\!+\!1)
\cdot e^{\scriptstyle 1\over\scriptstyle 2q}}\, 
$$
recalling $\nu=2\varpi/\sigma^2\!-\!1$. With $\rho$ any non-negative 
real, the {\it trigono\-me\-tric terms\/} $\ctrig {\xi} (\rho) $  
with $\xi$ equal to $\nu$ or $\nu\!+\!2$ are the integrals
$$
\ctrig {\xi} (\rho) = \int_0^{{\pi\over 2}}
\re\bigg( H_{-(\nu+4)}\bigg(
- {\cosh(\rho\!+\!i\phi)\over \sqrt{2q\,}}\bigg)
E_\xi(h)(\rho\!+\!i\phi)\bigg) d\phi
$$
over the real parts of products of Hermite functions times 
certain functions $E_b(h)$, and the {\it hyperbolic terms\/} 
$\chyp {\xi} (\rho) $ with $\xi$ equal to $\pm \nu$ or 
$\pm (\nu\!+\!2)$ are the integrals
$$
\chyp {\xi} (\rho) =
\int_\rho^\infty 
\im\bigg(H_{-(\nu+4)}\bigg(
- {\sinh(y)\over \sqrt{2q\,}}\, i\bigg)
E_\xi(h)(y\!+\!i{\pi\over 2})\bigg) dy
   $$
over the imaginary parts of such products. Herein $E_\xi(h)$ are the 
weighted complementary error functions for any complex $w$ given by
$$
E_\xi(h)(w)=
e\vbox to 9pt{}^{w\xi}\erfc\bigg(
   {w\over \sqrt{2h\,}}+ {\xi\over 2}\sqrt{2h\, }\bigg)\,  .
$$
\kleinerAbsatz
{\bf Remark:}\quad If $\rho$ equals zero, the trigonometric terms 
specialize to 
$$
\ctrig {\xi} (0) =
2\int_0^{\pi\over 2} H_{-(\nu+4)} \left( -{\cos (\phi)\over \sqrt{2q}}\right)
\cos(\xi\phi)\, d\phi\, . 
$$  
\mittlererAbsatz
Comparing with the formula of \S 5, the above formula,  
is given as a sum of single integrals whose integrands have a 
structural interpretation as products of two functions. It identifies 
the higher transcendental functions occuring as factors in these 
products, and  shows how they are given by or built up from Hermite 
functions. 
\goodbreak\grosserAbsatz
{\large\bf 24.\quad Epilogue:}\quad
On a technical level, the differences just noted between \S 5's and 
\S 23's formula can be regarded as consequences 
of the different mathematical approaches for proving the valuation 
formula. In fact,  \S 5's originates from Yor's direct attack on the 
law of the integral of geometric Brownian motion. In contrast, \S 23's 
is eventually based on the indirect enveloping construction of \S 7.
In fact, to obtain Asian option prices, first analytically invert 
\S 8~Theorem's Laplace transforms $F_{GY,a}$, then proceed using the 
key reduction of \S 7~Lemma;  formally speaking:
$
C^{(\nu)}(h,q)= \la^{-1}( F_{GY,q})(h)
$. 
\mittlererAbsatz
However, there are not only structural differences between \S5's 
and \S 23's formula. Finally deriving benchmarks for the normalized 
prices $C^{(\nu)}(h,q)$ of Asian options seems to be one of the 
main practical application of such formulas. 
As we have already mentioned, Yor's formula seems impracticable 
for this purpose in particular because of the gigantic size of the 
numbers it involves. For instance, consider valuing Asian 
options with annual interest rates $r$ equal to nine percent, 
with maturities of one year, with $K=S_0$, with $t_0=0$, and 
with a volatility of $\sigma=30\%$. These values require coping 
with numbers of order $10^{100}$. 
Our formula improves on this aspect as well. In 
the above situation, for instance, the four hyperbolic terms have orders
$10^7$ and the two trigonometric terms have orders $10^{-2}$. Sharpening
results of \cite{Rogers+Shi}, it thus became possible to derive in 
\cite[Chapter 5]{SHabil} the following first time benchmark values 
for normalized Asian option prices
$$ 
\vcenter{
\vbox{\offinterlineskip 
\hrule
\halign{&\vrule $#$& \strut \enspace  $\hfil#\hfil$\enspace\cr
&\sigma &&
{\rm maximal\, error}&&C^{\raise-1.5pt\hbox{$\scriptstyle (\nu)$}}(h,q)&\cr 
\noalign{\hrule}
&20\% && 4\,.\,9727 \times 10^{-16}
&&0\,.\,00074155998788343 &\cr 
\noalign{\hrule}
&30\% && 4\,.\,9687 \times 10^{-16}
&&0\,.\,00217354504625037&\cr
\noalign{\hrule}
&40\% &&4\,.\,9157 \times 10^{-16}
&&0\,.\,00478100328341654 &\cr 
\noalign{\hrule}
&50\% && 4\,.\,9461\times 10^{-16} 
&&0\,.\,00890942045213227 &\cr 
\noalign{\hrule}
}}}
$$
\centerline{\hbox{{\ninebf Table 3.}\enspace 
\ninerm Normalized prices $\scriptstyle C^{(\nu)}(h,q)$ of 
the Asian option for 
$\scriptstyle T=1$}.}
\mittlererAbsatz
All of this may be seen as a consequence of the Laplace transform 
approach to valuing non-Asian options initiated in \cite{Geman+Yor}. 
With hindsight, this paper appears to be a rich source for new results
in both mathematics and finance.
\grosserAbsatz
%
%
\vbox{\baselineskip=9truept\ninerm\hskip-23pt 
{\ninebf Acknowledgements:} The second author is grateful to the support
by the {\it Deutsche Forschungsgemeinschaft\/} and the hospitality of 
the {\it Mathematisches Institut\/} of the {\it Universit\"at Mannheim\/}.
While the usual disclaimer applies, we would also like to thank the 
editor Professor Shiryaev for his interest in the paper and the efforts 
he has spent on it.} 
%
%
\font\ninett=pcrb at 9pt 


\begin{thebibliography}{20}
%
\bibitem[\bf B]{Beals} R. Beals, {\nineit Advanced mathematical 
        analysis\/}, GTM 12, Springer 1973.
%
\bibitem[\bf CS]{carr-sch} P. Carr, M. Schr\"oder, On the 
        valuation of arithmetic-average Asian options: the Geman-Yor
         Laplace transform revisited, Mannheim and New York, December 2000,
        {\ninett http://arXiv.org/abs/math.CA/0102080}.
\bibitem[\bf C]{Conway} J.B. Conway, {\nineit Functions of one 
        complex variable\/}, 2nd ed., Springer 1984.
\bibitem[\bf D]{Doetsch} G. Doetsch, {\nineit Handbuch der Laplace 
       Transformation\/} I, Birkh\"auser Verlag 1971.
\bibitem[\bf DGY]{DGY} C. Donati-Martin, R. Ghomrasni, M. Yor:
         On certain Markov processes attached to exponential functionals 
         of Brownian motion: application to Asian options, {\nineit
         Revista Mathematica Iberoamericana\/} {\bf 17} (2001), 179--193.  
\bibitem[\bf D88]{Duffie88} D. Duffie: {\nineit Security markets\/},
         Academic Press 1988.
\bibitem[\bf D96]{Duffie92} D. Duffie: {\nineit Dynamic asset pricing 
        theory\/}, 2nd ed., Princeton UP 1996.
\bibitem[\bf Du90]{Dufresne90} D. Dufresne: The distribution of a 
        perpetuity: with applications to risk theory and pension funding, 
        {\nineit Scand. Actuarial J.\/} (1990), 39--79.
\bibitem[\bf Du00]{Dufresne00} D. Dufresne: Laguerre series for Asian and
        other options, {\it Math. Finance\/} {\ninebf 10} (2000), 407--28.  
\bibitem[\bf FB]{Freitag} E. Freitag, R. Busam, 
     {\nineit Funktionentheorie\/}, Springer 1993.
\bibitem[\bf FMW]{Fu+Madan+Wang} M.C. Fu, D.B. Madan, T. Wang: Pricing
     continuous Asian options: a comparison of Monte Carlo and Laplace 
     inversion methods, {\nineit J. Comp. Fin.\/} {\ninebf 2} (1998), 49--74.
\bibitem[\bf GY]{Geman+Yor} H. Geman, M. Yor: Bessel processes, Asian 
    options, and perpetuities, {\nineit Math. Finance\/} {\ninebf 3} (1993), 
    349--375.
\bibitem[\bf KS]{KSb} I. Karatzas, S.E. Shreve: {\nineit Methods
of mathematical finance\/}, Springer 1998.
\bibitem[\bf K]{Knight} F.B. Knight: {\nineit Essentials of Brownian
     motion and diffusion\/}, AMS Mathematical surveys 18, Providence 1991.
%
\bibitem[\bf L]{Lebedev}  N.N. Lebedev:  {\nineit Special functions and
       their applications\/}, Dover Publications 1972. 
\bibitem[\bf MR]{MR} M. Musiela, M. Rutkowski: {\nineit Martingale
        methods in financial modelling\/}, Springer 1997.
\bibitem[\bf \O]{Oeksendal} B. \O ksendal: {\nineit Stochastic 
        differential equations,} 5th. ed., Springer 1998.
\bibitem[\bf R]{Rudin} W. Rudin, {\nineit Real and complex analysis}, 
        3rd ed., McGraw Hill 1987.
\bibitem[\bf RY]{Revuz+Yor} D. Revuz, M. Yor: {\nineit Continuous martingales 
        and Brownian motion}, 2nd ed., Springer 1994.
\bibitem[\bf RS]{Rogers+Shi} L.C.G. Rogers, Z. Shi: The value of an Asian 
        option, {\nineit J. Appl. Probability\/} {\ninebf 32} (1995), 
        1077-1088.
\bibitem[\bf SA]{Asia} M. Schr\"oder: On the valuation of
    arithmetic-average Asian options: integral representations, 
    Oktober 1997, revised November 1999,
{\ninett http://arXiv.org/abs/math.CV/0003055}.
\bibitem[\bf SE]{Asia-Rechenbsp} M. Schr\"oder: On the valuation of
    arithmetic-average Asian options: explicit formulas, Universit\"at
    Mannheim, M\"arz 1999.
\bibitem[\bf SH]{SHabil} M. Schr\"oder: {\nineit Mathematical 
     ramifications of option valuation: the case of the Asian option\/},
     Habilitationsschrift Universit\"at Mannheim, April 2002. 
\bibitem[\bf SL]{Asia-Laguerre} M. Schr\"oder: On the valuation of
    arithmetic-average Asian options: Laguerre series and Theta 
    integrals, Dezember 2000,
      {\ninett http://arXiv.org/abs/math.CA/0012072}.
\bibitem[\bf W]{Watson} G.N. Watson: {\nineit A treatise on the theory of
    Bessel functions}, 2nd ed., Cambridge University Press 1944. 
\bibitem[\bf We]{Weil} A. Weil: {\nineit \OE uvres scientifiques\/}, 
    volume 3, Springer. 
\bibitem[\bf Y80]{Yor 80}  M. Yor: Loi d'indice du lacet Brownien, et
    distribution de Hartman-Watson, 
    {\nineit Z. Wahr\-schein\-lichkeitstheorie\/} {\ninebf 53} (1980), 
    71--95.
\bibitem[\bf Y92]{Yor 92a}  M. Yor: Sur certaines fonctionnelles 
     exponentielles du mouvement Brownien r\'eel, {\nineit J. Appl. Prob.\/} 
     {\ninebf 29} (1992), 202--208. 
\bibitem[\bf Y]{Yor 92}  M. Yor: On some exponential functionals of 
    Brownian motion, {\nineit Adv. Appl. Prob.} {\ninebf 24} (1992), 
    509--531.
\bibitem[\bf YG\"o]{Goeing-Yor}: M. Yor, A. G\"oing-Jaeschke: A survey
    and some generalizations of Bessel processes,  ETH Z\"urich 1999.
\bibitem[\bf Y00]{Yor 3}  M. Yor et al.: {\nineit Exponential 
    functionals of Brownian motion and related processes III\/}, 
    pre-print Paris VI, May 2000.
\bibitem[\bf Y01]{Yor01}  M. Yor et al.: {\nineit Exponential 
    functionals of Brownian motion and related processes\/}, 
    Springer 2001.
\end{thebibliography}
\end{document}